\documentclass[conference]{new-aiaa}
\usepackage[utf8]{inputenc}
\usepackage{textcomp}

\usepackage{graphicx}
\usepackage{subcaption}
\DeclareCaptionLabelFormat{aiaa}{#2)}
\usepackage{amsmath}
\usepackage[version=4]{mhchem}
\usepackage{siunitx}
\usepackage{longtable,tabularx}
\usepackage[printonlyused]{acronym}
\acrodef{uut}[UUT]{Unified Unscented Transform}
\acrodef{ut}[UT]{Unscented Transform}
\acrodef{mc}[MC]{Monte Carlo}
\acrodef{bvp}[BVP]{Boundary Value Problem}
\acrodef{gut}[GUT]{Generalized Unscented Transform}
\acrodef{cut}[CUT]{Conjugate Unscented Transform}
\acrodef{house}[HOUSE]{Higher-Order Unscented Estimator}
\acrodef{mce}[MCE]{Moment Constraint Equation}
\acrodef{ug}[UG]{Unscented Guidance}
\acrodef{dof}[DoF]{Degrees of Freedom}
\acrodef{stargate}[STARGATE]{STochastic Algorithm for Robust Guidance Analysis and Trajectory Estimation}
\acrodef{nlp}[NLP]{Nonlinear Programming}
\acrodef{raptor}[RAPTOR]{Rapid Pseudospectral Trajectory Optimization and Robustification}
\acrodef{ocp}[OCP]{Optimal Control Problem}
\acrodef{lgr}[LGR]{Legendre-Gauss-Radau}
\acrodef{lgl}[LGL]{Legendre-Gauss-Lobatto}
\acrodef{flgr}[fLGR]{flipped-Legendre-Gauss-Radau}
\acrodef{rlv}[RLV]{Reusable Launch Vehicle}
\acrodef{ps}[PS]{Pseudospectral}
\acrodef{ipopt}[IPOPT]{Interior Point OPTimizer}
\acrodef{snopt}[SNOPT]{Sparse Nonlinear OPTimizer}
\acrodef{worhp}[WORHP]{We Optimize Really Huge Problems}
\acrodef{rtls}[RTLS]{Return To Launch Site}
\acrodef{kkt}[KKT]{Karush-Kuhn-Tucker}
\acrodef{erc}[ERC]{European Research Council}
\acrodef{pjrn}[PJRN]{Projected Jacobian Row Normalization}
\acrodef{cut4}[CUT4]{Fourth-Order Conjugate Unscented Transform}
\hypersetup{hidelinks}
\usepackage{dsfont}
\usepackage{todonotes}
\usepackage[table]{xcolor}
\usepackage{algorithm}
\usepackage{algorithmic}

\title{A Unified Unscented Transform for Uncertainty Quantification and Propagation}

\author{Riccardo Minnozzi\footnote{PhD Student, Department of Industrial Engineering, Via Luciano Montaspro 97, \href{mailto:riccardo.minnozzi@unibo.it}{riccardo.minnozzi@unibo.it}.} and Marco Sagliano \footnote{Full Professor, Department of Industrial Engineering, Via Luciano Montaspro 97, \href{mailto:marco.sagliano@unibo.it}{marco.sagliano@unibo.it}, AIAA Associate Fellow}}
\affil{University of Bologna, Forlí, 47121, Italy}

\begin{document}

\maketitle
\begin{abstract}
Unscented Optimal Control replaces expectation integrals in uncertain systems with sigma-point rules, whose points and weights set both the size of the problem and the accuracy of its uncertainty model. Existing rules adopt fixed-point structures, which may cause loss of accuracy in generic non-Gaussian distributions, violation of support-set constraints, or simply introduce more sigma-points than strictly required. This paper presents the Unified Unscented Transform: the sigma-point rule is designed through a constrained optimization problem in the point locations and associated weights. We enforce the moment-matching for a generic distribution while taking into account  additional validity constraints. We then propose a sparsification mechanism that pushes for sparse solutions of the same problem without significant loss of accuracy. We perform a numerical assessment of the technique over a broad range of multidimensional distributions, showing that it outperforms existing sigma-point construction approaches in terms of both sparsity and accuracy in most of the tested cases. Finally, we formulate an entry problem under uncertain atmospheric density within the Unscented Optimal Control framework to show the ability of our methodology to accurately handle non-Gaussian distributions and chance constraints.
\end{abstract}

\section{Introduction}
\label{sec:introduction}

\par Uncertainty is an inherent feature of space missions, and treating it just as a disturbance to be rejected might not only lead to conservative trajectory designs, but also places high demands on feedback controllers when significant disturbances occur \cite{minnozzi2026_STARGATERobustGuidance}. This motivates robust and stochastic trajectory-planning frameworks that include uncertainty directly in the decision-making process, such that performance, robustness, and risk can be balanced before execution. 

\par Stochastic optimal control, which accounts explicitly for random quantities in the dynamics, objectives, or constraints \cite{astrom_2006_IntroductionStochastic,bertsekas_2007_StochasticOptimal,fleming_1975_DeterministicStochastic}, has been applied to trajectory design through augmented shadow trajectories \cite{grant2015_RapidRobustTrajectory}, linear covariance propagation \cite{ridderhof2020_ChanceconstrainedCovarianceControl, benedikter2022_ConvexApproachCovariance}, Desensitized Optimal Control \cite{seywald2019_DesensitizedOptimalControl, robbiani2024_FastDesensitizedOptimal}, and distributional approximations such as Gaussian mixture models \cite{vittaldev2016_SpacecraftUncertaintyPropagation}. Each formulation balances the fidelity of the uncertainty model against the size and conditioning of the resulting optimization problem. For nonlinear Stochastic Optimal Control, three families of methods are particularly relevant to the present work: sampling, polynomial representations, and unscented methods.

\par Sampling-based methods such as \ac{mc} estimate the statistics of interest from a set of realizations and are the natural reference for validation, but the number of samples required for chance constraints and tail probabilities makes them unfit for repeated or onboard optimization \cite{telen2015_ApproximateRobustOptimization, ross2015_RiemannStieltjesOptimal}. Polynomial-chaos methods replace the stochastic variables with orthogonal polynomial expansions \cite{xiu_2010_NumericalMethods} and have been applied to trajectory generation under uncertainty \cite{fisher2011_OptimalTrajectoryGeneration, whittle2018_StochasticOptimalTrajectory}, but their cost grows quickly with the expansion order and the number of uncertain variables, which adaptive constructions only alleviate \cite{jiang_2020_HighdimensionalUncertainty}.

\par Unscented methods approximate the input distribution through a weighted set of sigma points, propagate each one through the nonlinear model, and reconstruct the output statistics from the propagated set \cite{julier2004_UnscentedFilteringNonlinear}. The approach is non-intrusive as in \ac{mc} sampling, but replaces the large sample set with a sigma-point rule: a small, carefully constructed set of points and weights that acts as a multivariate quadrature rule (i.e., a cubature rule) for the input distribution. Unscented Optimal Control follows this strategy by defining a common control policy over each point of this quadrature rule, such that the sigma points set both the size and accuracy of the problem \cite{ross2014_UnscentedOptimalControl, ross2015_RiemannStieltjesOptimal, ross2025_UnscentedTrajectoryOptimization}. The compactness and non-intrusiveness of this representation motivate the focus of this paper on unscented methods.

\par Multiple strategies to construct the unscented rule are available in the literature, and we focus on a representative subset in this work. For $n$ uncertain variables, the \ac{ut} constructs $2n+1$ sigma points that reproduce the mean and covariance of the input distribution \cite{julier2004_UnscentedFilteringNonlinear}, and it has been applied extensively to navigation and guidance \cite{wan2000_UnscentedKalmanFilter,schwartz2025_UnscentedHypersonicTrajectory}. The \ac{cut} places additional sigma points along conjugate axes to reproduce higher-order moments of Gaussian and uniform distributions \cite{adurthi2012_ConjugateUnscentedTransform, adurthi2013_ConjugateUnscentedTransform}, at the cost of a number of points that grows exponentially with $n$. It has also been applied to chance-constrained optimal control \cite{burnett2026_AnalyticNonGaussianConfidence}. The \ac{gut} retains the $2n+1$ points of the \ac{ut}, with asymmetric offsets that match only the diagonal entries of the skewness and kurtosis tensors \cite{ebeigbe2021_GeneralizedUnscentedTransformation}. Each construction therefore fixes the matched moment entries, the number of points, and the admissible distributions, such that a rule may miss the moments that matter for a non-Gaussian distribution, or use more points than these moments require.

\par The resulting rule, moreover, is not always valid (i.e., with nonnegative weights and every point inside the support of the input distribution): the central weight of the \ac{ut} \cite{julier2000_NewMethodNonlinear} and of the \ac{gut} (Sec.~\ref{subsec:existing-rules}) can become negative, and the points of the \ac{ut} and the \ac{cut} can leave the bounded support of the specific distribution of interest. Both limitations matter in the \ac{ug} framework of Refs.~\cite{ross2015_RiemannStieltjesOptimal, ross2025_UnscentedTrajectoryOptimization}, where the weights must form a probability measure for the expectation integrals, and each point outside the support adds a physically meaningless trajectory. Instead of prescribing a sigma-point structure independently of the target distribution and computational budget, the \ac{uut}\footnote{available at: \url{https://github.com/RiccardoMinnozzi-UniBo/Unified-Unscented-Transform}} introduced in this paper formulates the rule design as a constrained optimization problem: as in designed quadrature \cite{keshavarzzadeh2018_NumericalIntegrationMultiple}, the locations and weights of the sigma points are the decision variables, while the validity conditions and the prescribed moments are the constraints. A sparsification procedure then reduces the number of points that a valid and accurate rule requires.

\par The principal contributions of this paper are the following. First, we introduce the \ac{uut} numerical cubature framework, consisting of the constrained optimization of sigma points and weights to match prescribed moment entries. Second, we propose a sparsification mechanism, combining reweighted $L_1$ minimization and greedy point elimination to reduce the number of points required for a valid rule. Third, we carry out a systematic comparison against existing unscented methods (\ac{ut}, \ac{cut}, and \ac{gut}), showing that the existing rules are particular solutions of the \ac{uut} problem wherever they are verified, and that the \ac{uut} can match the same moment conditions with fewer points, or more moment conditions with the same number of points. Finally, we integrate the \ac{uut} into the \ac{ug} framework and showcase the effect of higher-order statistical information on quantile approximations for the chance-constrained problem. The push for efficient robust optimal control methods is part of the broader \acs{erc}-funded \acs{stargate} project \cite{minnozzi2026_STARGATERobustGuidance}.

\par The remainder of the paper is organized as follows. Section~\ref{sec:unified-unscented-transform} formulates the \ac{uut} rule design problem, starting from the moment-matching conditions satisfied by existing rules. Section~\ref{sec:rule-sparsification} presents the sparsification mechanism. Section~\ref{sec:numerical-assessment} compares the \ac{uut} against existing unscented methods across multiple distributions and dimensions, and evaluates the effectiveness of the sparsification technique. Section~\ref{sec:unscented-guidance} integrates the \ac{uut} into the \ac{ug} framework and applies it to an entry problem under atmospheric density uncertainty. Finally, Section~\ref{sec:conclusions} summarizes the main findings of this work. 
\section{The Unified Unscented Transform}
\label{sec:unified-unscented-transform}

\par As discussed in Sec.~\ref{sec:introduction}, existing sigma-point rules share a common approach: each rule assumes a specific point structure (e.g., symmetric points along the principal axes for the \ac{ut}, or additional points along conjugate axes for the \ac{cut}) and then determines the point locations and weights that match a set of prescribed statistical moments. The \ac{uut} retains the moment conditions and removes the structural assumption: point locations and weights are both treated as decision variables of a constrained optimization problem, whose constraints include the prescribed moments and the validity of the rule. This section presents the \ac{uut} technique, starting from the moment-matching conditions and their relation to existing rules, then moving onto the rule design problem, and finally addressing its scaling.

\subsection{Sigma-Point Rules as Moment-Matching Cubature}
\label{subsec:moment-matching-cubature}

\par Consider the vector of uncertain variables $\mathbf{x} \in \mathbb{R}^n$, distributed according to the target distribution, which the rule approximates. We assume its support to be contained in the box $[\mathbf{a}, \mathbf{b}]$, such that $\operatorname{supp}(\mathbf{x}) \subseteq [\mathbf{a}, \mathbf{b}]$, where the entries of $\mathbf{a}$ and $\mathbf{b}$ may be infinite (e.g., for a Gaussian distribution). A sigma-point rule of $N$ points approximates the expectation of a function $g(\mathbf{x})$ through the weighted sum
\begin{equation}
    \mathbb{E}\left[g(\mathbf{x})\right] \approx \sum_{i=1}^{N} w_i\, g(\mathbf{x}_i)
    \label{eq:sigma-point-rule}
\end{equation}
where $\mathbf{x}_i \in \mathbb{R}^n$ and $w_i \in \mathbb{R}$ are the $i$-th sigma point and its weight, collected in the matrix $\mathbf{X} = [\mathbf{x}_1, \dots, \mathbf{x}_N] \in \mathbb{R}^{n \times N}$ and in the vector $\mathbf{w} = [w_1, \dots, w_N]^T \in \mathbb{R}^N$, respectively. In unscented methods, $g(\mathbf{x})$ represents the nonlinear model through which each sigma point is propagated, and Eq.~\eqref{eq:sigma-point-rule} is the cubature rule that reconstructs the statistics of its output. For a sigma-point rule to represent the target distribution, its weights must be nonnegative and sum to one (i.e., unit mass), and its points must lie inside the support, therefore
\begin{equation}
    w_i \geq 0, \qquad \sum_{i=1}^{N} w_i = 1, \qquad \mathbf{a} \leq \mathbf{x}_i \leq \mathbf{b}, \qquad i = 1, \dots, N
    \label{eq:valid-rule}
\end{equation}
where the inequalities on the points hold componentwise, as for the support box. We refer to a rule that satisfies these conditions as valid: within an unscented method, nonnegative weights prevent negative variance estimates, and points inside the support keep every propagated point physically meaningful.

\par The statistical information of the target distribution is expressed through its statistical moments: the mean $\boldsymbol{\mu} = \mathbb{E}[\mathbf{x}]$ and the central moments of order $k \geq 2$, whose entries are defined as
\begin{equation}
    m_{\mathbf{j}} = \mathbb{E}\left[\prod_{l=1}^{k} \left(x_{j_l} - \mu_{j_l}\right)\right], \qquad \mathbf{j} = (j_1, \dots, j_k)
    \label{eq:target-moments}
\end{equation}
where the scalar index $j_l$, the $l$-th entry of the tuple $\mathbf{j}$, selects the component of $\mathbf{x}$ in the $l$-th of the $k$ factors, and the length $|\mathbf{j}| = k$ is the order of the entry. The entries in Eq.~\eqref{eq:target-moments} form the central moment tensors $\mathcal{M}_k$: $\mathcal{M}_2 = \mathbf{P}$ is the covariance matrix, while $\mathcal{M}_3$ and $\mathcal{M}_4$ are the third- and fourth-order tensors, referred to as the skewness and kurtosis tensors. Each entry $m_{\mathbf{j}}$ occupies the position $(j_1, \dots, j_k)$ of $\mathcal{M}_k$: for instance, $m_{(1,2)} = [\mathcal{M}_2]_{12} = P_{12}$ is the covariance between $x_1$ and $x_2$. For uniformity of notation, we also set $\mathcal{M}_1 = \boldsymbol{\mu}$, with entries $m_{(j)} = \mu_j$. The moments of a sigma-point rule are defined analogously by replacing expectations with the weighted sum of Eq.~\eqref{eq:sigma-point-rule}: the rule mean is $\hat{\boldsymbol{\mu}} = \mathbf{X}\mathbf{w}$, and the entries of the rule tensors $\hat{\mathcal{M}}_k$ are
\begin{equation}
    \hat{m}_{\mathbf{j}}(\mathbf{X}, \mathbf{w}) = \sum_{i=1}^{N} w_i \prod_{l=1}^{k} \left(x_{i,j_l} - \hat{\mu}_{j_l}\right)
    \label{eq:rule-moments}
\end{equation}
where $x_{i,j}$ denotes the $j$-th component of the $i$-th sigma point $\mathbf{x}_i$.

\par Since the moment tensors are symmetric, their entries $m_{\mathbf{j}}$ do not change when the indices in $\mathbf{j}$ are permuted. For the \ac{uut} framework we avoid redundant moment conditions by only considering the index tuples with $j_1 \leq j_2 \leq \dots \leq j_k$, such that each moment tensor of order $k$ contributes $\binom{n+k-1}{k}$ unique entries. The moment-matching conditions are imposed on a prescribed set $\mathcal{J}$ of unique entries of order up to $K$ through
\begin{equation}
    \hat{m}_{\mathbf{j}}(\mathbf{X}, \mathbf{w}) = m_{\mathbf{j}}, \qquad \mathbf{j} \in \mathcal{J}
    \label{eq:moment-conditions}
\end{equation}
Two prescribed sets are used throughout this paper. The full set $\mathcal{J}_K^{\mathrm{f}}$ contains every unique entry up to order $K$. The marginal set $\mathcal{J}_K^{\mathrm{m}}$ contains every entry of orders one and two, but only the marginal entries $\mathbf{j} = (j, \dots, j)$ of the higher orders, meaning that the mixed entries involving more than one variable are ignored above the covariance. Table~\ref{tab:moment-entries} provides an example of the indexing conventions for a distribution with $n=2$ components and moments up to order $K=3$: the tuples $\mathbf{j}$ of length one select the components of the mean, those of length two the entries of the covariance matrix, and those of length three the third-order entries, which belong to the marginal set only when all their indices are equal. 

\par Matching the full set $\mathcal{J}_K^{\mathrm{f}}$ together with the unit mass is equivalent to integrating exactly, through Eq.~\eqref{eq:sigma-point-rule}, every polynomial in $\mathbf{x}$ of degree up to $K$, since the expectation of each monomial follows from the mean and the central entries up to order $K$. A rule for a full set is therefore a positive cubature rule of degree $K$, to which the bounds on the number of points of Secs.~\ref{subsec:reweighted-l1} and~\ref{subsec:greedy-elimination} apply.

\begin{table}[htb]
\caption{\label{tab:moment-entries} Unique moment entries for $n = 2$ and $K = 3$, and their membership in the full and marginal sets.}
\centering
\begin{tabular}{llllcc}
\hline\hline
Order $k$ & Tuple $\mathbf{j}$ & Entry $m_{\mathbf{j}}$ & Type & $\mathcal{J}_3^{\mathrm{f}}$ & $\mathcal{J}_3^{\mathrm{m}}$ \\
\hline
1 & $(1)$, $(2)$ & $\mu_1$, $\mu_2$ & Mean & Yes & Yes \\
2 & $(1,1)$, $(2,2)$ & $P_{11}$, $P_{22}$ & Variance & Yes & Yes \\
2 & $(1,2)$ & $P_{12}$ & Covariance & Yes & Yes \\
3 & $(1,1,1)$ & $\mathbb{E}[(x_1-\mu_1)^3]$ & Marginal & Yes & Yes \\
3 & $(2,2,2)$ & $\mathbb{E}[(x_2-\mu_2)^3]$ & Marginal & Yes & Yes \\
3 & $(1,1,2)$ & $\mathbb{E}[(x_1-\mu_1)^2(x_2-\mu_2)]$ & Mixed & Yes & No \\
3 & $(1,2,2)$ & $\mathbb{E}[(x_1-\mu_1)(x_2-\mu_2)^2]$ & Mixed & Yes & No \\
\hline\hline
\end{tabular}
\end{table}

\par Table~\ref{tab:existing-rules} summarizes how the existing rules fit this framework. The \ac{ut} places $2n$ points symmetrically along the columns of a matrix square root of $\mathbf{P}$ and one point on the mean, a construction that reproduces the mean and covariance of any distribution ($\mathcal{J}_2^{\mathrm{f}}$) and, by symmetry, the vanishing third-order moments of symmetric distributions ($\mathcal{J}_3^{\mathrm{f}}$) \cite{julier2004_UnscentedFilteringNonlinear}. The \ac{gut} retains the same $2n+1$ points but places them asymmetrically, such that the marginal third- and fourth-order entries are also matched ($\mathcal{J}_4^{\mathrm{m}}$) \cite{ebeigbe2021_GeneralizedUnscentedTransformation}. The \ac{cut} adds points along conjugate axes, and the \ac{cut4} reproduces every entry up to order four (and also order five, by symmetry) for Gaussian distributions ($\mathcal{J}_5^{\mathrm{f}}$) at the cost of $2n + 2^n$ points for $n \geq 3$ \cite{adurthi2012_ConjugateUnscentedTransform}. A second \ac{cut4} construction, with its own radii, weights, and conjugate axes, matches the same set for independent uniform variables, with tabulated rules for 2 to 8 variables \cite{adurthi2013_ConjugateUnscentedTransform}. The higher-order rules of the \ac{cut} extend the same construction to the higher-order full sets for Gaussian distributions, at the price of a significantly larger number of points. Given its use in chance-constrained optimal control \cite{burnett2026_AnalyticNonGaussianConfidence}, this paper adopts the \ac{cut4} as the reference higher-order sigma-point rule.

\par As shown by Table~\ref{tab:existing-rules}, the number of points of these rules follows from their structure rather than their prescribed set, and the two numbers grow at different rates: the \ac{ut} and \ac{gut} use a number of points that grows linearly with $n$ for a number of moment conditions that grows cubically and quadratically, respectively, while the points of the \ac{cut4} grow exponentially with $n$ for conditions that grow as $n^5$. The \ac{uut} instead solves Eq.~\eqref{eq:moment-conditions} directly for the point locations and weights, such that the number of points can be chosen independently of the prescribed set. Consistently with this interpretation, Sec.~\ref{subsec:existing-rules} shows that the existing rules are particular solutions of the \ac{uut} problem wherever they are verified (i.e., valid and accurate on their prescribed set), while Sec.~\ref{sec:rule-sparsification} discusses how few points are required for a specific prescribed set.

\begin{table}[htb]
\caption{\label{tab:existing-rules} Existing sigma-point rules as particular solutions of the moment-matching conditions, with their numbers of points and conditions for $n$ variables.}
\centering
\renewcommand{\arraystretch}{1.3}
\begin{tabularx}{\textwidth}{l>{\raggedright\arraybackslash}X>{\raggedright\arraybackslash}p{3.4cm}cc}
\hline\hline
Rule & Point structure & Prescribed set & Points $N$ & Conditions $|\mathcal{J}|+1$ \\
\hline
\acs{ut} & Mean and symmetric pairs along the columns of $\mathbf{P}^{1/2}$ & $\mathcal{J}_3^{\mathrm{f}}$ (symmetric targets) & $2n+1$ & $\binom{n+3}{3}$ \\
\acs{gut} & Mean and asymmetric pairs along the columns of $\mathbf{P}^{1/2}$ & $\mathcal{J}_4^{\mathrm{m}}$ (any target) & $2n+1$ & $(n^2+7n+2)/2$ \\
\acs{cut4}, Gaussian & Symmetric points along principal and conjugate axes & $\mathcal{J}_5^{\mathrm{f}}$ (Gaussian targets) & $2n+2^n$\textsuperscript{a} & $\binom{n+5}{5}$ \\
\acs{cut4}, uniform & Symmetric points along principal and conjugate axes, with other radii and weights & $\mathcal{J}_5^{\mathrm{f}}$ (independent uniform targets) & Tabulated\textsuperscript{b} & $\binom{n+5}{5}$ \\
\hline\hline
\multicolumn{5}{l}{\textsuperscript{a} For $n \geq 3$ (9 points for $n = 2$).} \\
\multicolumn{5}{l}{\textsuperscript{b} For $2 \leq n \leq 8$: 8, 14, 24, 252, and 1808 points for $n = 2$, 3, 4, 6, and 8.}
\end{tabularx}
\end{table}

\subsection{Rule Design Problem}
\label{subsec:rule-design-problem}

\par The \ac{uut} computes a rule of $N$ points by solving the  \ac{nlp} problem
\begin{equation}
    \begin{aligned}
        \min_{\mathbf{X}, \mathbf{w}} \quad & J(\mathbf{X}, \mathbf{w}) \\
        \text{s.t.} \quad & \left|\hat{m}_{\mathbf{j}}(\mathbf{X}, \mathbf{w}) - m_{\mathbf{j}}\right| \leq \varepsilon_m, & & \mathbf{j} \in \mathcal{J} \\
        & \sum_{i=1}^{N} w_i = 1 \\
        & 0 \leq w_i \leq 1, \qquad \mathbf{a} \leq \mathbf{x}_i \leq \mathbf{b}, & & i = 1, \dots, N
    \end{aligned}
    \label{eq:rule-design-problem}
\end{equation}
where the point locations and weights are the decision variables. Non-trivial objectives are introduced in Sec.~\ref{sec:rule-sparsification}, while $J=0$ makes Eq.~\eqref{eq:rule-design-problem} a pure feasibility problem. The moment conditions of Eq.~\eqref{eq:moment-conditions} are imposed within the $\varepsilon_m$ band instead of equalities, such that the \ac{nlp} solver does not stop prematurely when the conditions outnumber the unknowns, as rules that satisfy every condition can exist in that case as well (e.g., symmetric rules, whose odd-order conditions hold by construction). The remaining constraints enforce the validity conditions of Eq.~\eqref{eq:valid-rule} on every feasible rule: the support enters as simple bounds on the point coordinates, while the weight bounds and the unit mass make every feasible rule a positive cubature rule.

\par The number of points $N$, referred to as the budget, is an input of Eq.~\eqref{eq:rule-design-problem}, whose $N(n+1)$ unknowns (the coordinates and the weight of each point) can exceed the $|\mathcal{J}|+1$ conditions: in this case, the problem admits many solutions. We initialize the \ac{uut} problem with $N$ equally weighted (i.e., $w_i=1/N$) points, deterministically sampled from a Gaussian surrogate of the target distribution, as detailed in Appendix~\ref{subsec:initial-guess}. The problem in Eq.~\eqref{eq:rule-design-problem} is solved with \ac{ipopt} \cite{wachter2006_ImplementationInteriorpointFilter}, using the analytic derivatives of the moment conditions derived in Appendix~\ref{subsec:analytic-derivatives} and the settings reported in Sec.~\ref{subsec:experimental-setup}. Every sigma-point rule is checked a posteriori in the original variables: we refer to a rule as accurate when it matches every entry of its prescribed set within the verification tolerance $\varepsilon_v$, and as verified when it is both valid and accurate. Solving Eq.~\eqref{eq:rule-design-problem} costs more than evaluating the closed-form expressions of the existing rules, but this cost is paid once and offline: the rule depends only on the uncertainty model, and the guidance problem uses its points and weights as known data. The \ac{uut} construction times for all the studies and distributions considered in Sec.~\ref{sec:numerical-assessment} ranges from seconds to few minutes for the high-dimensional and highly constrained cases, while the iterative sparsification procedure requires longer runtimes (up to about 2 hours).

\subsection{Problem Scaling}
\label{subsec:problem-scaling}

\par The moment entries prescribed in Eq.~\eqref{eq:rule-design-problem} can differ by several orders of magnitude, and this difference degrades the numerical conditioning of the problem. An entry of order $k$ scales with the product of the $k$ standard deviations it involves: for instance, the fourth-order central moment of a Gaussian variable with standard deviation $\sigma$ is $3\sigma^4$. In guidance problems, the uncertain parameters may have different physical units, and two variables with $\sigma = 0.1$ and $\sigma = 10$ already produce fourth-order entries that are eight orders of magnitude apart within the same problem. Strong correlations add a second difficulty: as $\mathbf{P}$ approaches a singular matrix, the correlated components take similar values on every sigma point, such that the moment conditions of the same order (e.g., on $P_{11}$, $P_{12}$, and $P_{22}$) become nearly redundant and the problem becomes ill conditioned. To address these conditioning aspects, we solve Eq.~\eqref{eq:rule-design-problem} in \emph{whitened coordinates} (i.e., the variables for an equivalent distribution with zero mean and identity covariance). Introducing the Cholesky factorization $\mathbf{P} = \mathbf{L}\mathbf{L}^T$ of the target covariance, the whitened variables
\begin{equation}
    \mathbf{z} = \mathbf{L}^{-1}\left(\mathbf{x} - \boldsymbol{\mu}\right)
    \label{eq:whitening}
\end{equation}
have zero mean and identity covariance. Since Eq.~\eqref{eq:whitening} is linear, each central moment entry of the whitened distribution $\mathbf{z}$, denoted by the superscript $z$, follows from the entries of the same order of the target distribution $\mathbf{x}$ as
\begin{equation}
    m^{z}_{\mathbf{j}} = \sum_{\mathbf{q}} \left(\prod_{l=1}^{k} \left[\mathbf{L}^{-1}\right]_{j_l q_l}\right) m_{\mathbf{q}}
    \label{eq:whitened-moments}
\end{equation}
where the sum runs over all the index tuples $\mathbf{q}$ of length $k$, including the permutations of the unique ones. Each whitened entry is therefore a weighted sum of all the entries of the same order of $\mathbf{x}$, with one factor of $\mathbf{L}^{-1}$ for each of its $k$ indices. For $k = 2$, Eq.~\eqref{eq:whitened-moments} reduces to the familiar $\mathbf{P}^{z} = \mathbf{L}^{-1}\mathbf{P}\mathbf{L}^{-T} = \mathbf{I}$. Equation~\eqref{eq:rule-design-problem} is thus solved in the whitened coordinates $\mathbf{z}_i$, and the solution is mapped back through $\mathbf{x}_i = \mathbf{L}\mathbf{z}_i + \boldsymbol{\mu}$. Since Eq.~\eqref{eq:whitened-moments} is invertible order by order, a full set $\mathcal{J}_K^{\mathrm{f}}$ defines the same rules in both variables, while a marginal set is preserved only for independent variables, where $\mathbf{L}$ is diagonal and the whitening reduces to a standardization (i.e., individually mapping each component to a variable with zero mean and unit variance). Since a non-diagonal $\mathbf{L}^{-1}$ combines marginal and mixed entries of $\mathbf{x}$, marginal sets $\mathcal{J}_K^{\mathrm{m}}$ are solved in the original coordinates. For a nearly singular $\mathbf{P}$, the Cholesky factor $\mathbf{L}$ is computed accurately, and the large entries in its inverse $\mathbf{L}^{-1}$ amplify the rounding errors of the whitened targets of Eq.~\eqref{eq:whitened-moments}. While this amplification can theoretically make the whitened targets inconsistent with the original target distribution and cause solver failures, it is never verified in practice. Additionally, any identified rule is mapped back through $\mathbf{L}$ to the original variables, and the independent verification of Sec.~\ref{subsec:rule-design-problem} would filter out any accepted rule that violates the prescribed moment set (note that no such case was observed in practice throughout all the analyzed distributions). An exactly singular $\mathbf{P}$, instead, means that some components are linear combinations of the others, which shall be removed from $\mathbf{x}$ before the rule is constructed.

\par The support box of the original variables maps to a hyper-parallelepiped in the whitened coordinates, which cannot be represented by simple bounds on $\mathbf{z}_i$. The whitened variables are therefore bounded by the smallest box that encloses the mapped hyper-parallelepiped, while the original support is enforced exactly by the linear constraints
\begin{equation}
    \mathbf{a} \leq \mathbf{L}\mathbf{z}_i + \boldsymbol{\mu} \leq \mathbf{b}, \qquad i = 1, \dots, N
    \label{eq:whitened-support}
\end{equation}
which are added to Eq.~\eqref{eq:rule-design-problem}. Whitening improves the effectiveness of the solver on correlated distributions, without affecting the prescribed moments and support.
\section{Rule Sparsification}
\label{sec:rule-sparsification}

\par The number of sigma points sets the computational cost of every problem in which a rule is used. In \ac{ug} each sigma point propagates its own state trajectory through the dynamics, linearly increasing the problem size with $N$ \cite{ross2014_UnscentedOptimalControl}, and the rule design problem of Eq.~\eqref{eq:rule-design-problem} grows with $N$ as well. Since $N$ is a prescribed input of Eq.~\eqref{eq:rule-design-problem}, rules computed with overparameterized budgets can carry more points than required. While ideally sparsification would return the sparsest valid rule (i.e., the minimum number of active weights $w_i>0$), this is a combinatorial problem \cite{natarajan1995_SparseApproximateSolutions}, which is not solved exactly in this paper. We propose a two-stage procedure that reduces an overparameterized rule to a small verified one, without guaranteeing that smaller rules do not exist.

\subsection{Reweighted \texorpdfstring{$L_1$}{L1} Minimization with Fixed Points}
\label{subsec:reweighted-l1}

\par The standard approach to promote sparse solutions defines the objective of Eq.~\eqref{eq:rule-design-problem} as a regularization term, typically
\begin{equation}
    J_1(\mathbf{w}) = \lambda \sum_{i=1}^{N} |w_i| \qquad \text{or} \qquad J_2(\mathbf{w}) = \lambda \sum_{i=1}^{N} w_i^2
    \label{eq:plain-regularization}
\end{equation}
where $\lambda > 0$ is the regularization weight. The $L_1$ term in $J_1$ is the convex surrogate of the number of nonzero entries used in LASSO regression \cite{tibshirani2011_RegressionShrinkageSelection}, while the $L_2$ term in $J_2$ penalizes large entries and therefore favors a more uniform (although less sparse) solution. Neither term sparsifies a sigma-point rule: since the weights of a valid rule are nonnegative and sum to one, $J_1 = \lambda$ for every feasible rule and has no effect on the solution, while $J_2$ is minimized by equal weights $w_i = 1/N$ and therefore drives the rule away from sparsity, as demonstrated in Appendix~\ref{subsec:regularization-evidence}.

\par Reweighted $L_1$ minimization \cite{candes2008_EnhancingSparsityReweighted} replaces the constant penalty $J_1$ with a sequence of $L_1$ problems whose costs
\begin{equation}
    c_i^{(\ell)} = \frac{1}{w_i^{(\ell)} + 1/N}, \qquad i = 1, \dots, N
    \label{eq:reweighting-costs}
\end{equation}
are computed from the weights $w_i^{(\ell)}$ of the previous iteration $\ell$, with the offset $1/N$ keeping the costs bounded. Unlike $J_1$, the weighted sum $\sum_i c_i^{(\ell)} w_i$ is not constant on the feasible set: small weights receive high costs and are driven toward zero, while large weights remain almost free to change. Each iteration $\ell$ can be interpreted as a local linearization of the log-sum penalty $\sum_i \log(w_i + 1/N)$, a closer surrogate of the number of active weights than $J_1$ \cite{candes2008_EnhancingSparsityReweighted}.
\par The points are kept fixed during the reweighting, such that each iteration solves the linear program
\begin{equation}
    \begin{aligned}
        \min_{\mathbf{w}} \quad & \sum_{i=1}^{N} c_i^{(\ell)} w_i \\
        \text{s.t.} \quad & \mathbf{A}\mathbf{w} = \mathbf{A}\mathbf{w}^{(0)} \\
        & 0 \leq w_i \leq 1, \qquad i = 1, \dots, N
    \end{aligned}
    \label{eq:reweighting-lp}
\end{equation}
where $\mathbf{w}^{(0)}$ are the weights of the starting rule, and each row of the matrix
\begin{equation}
    \mathbf{A} = \begin{bmatrix} \mathbf{1}^T \\ \mathbf{X} \\ \mathbf{E} \end{bmatrix}, \qquad E_{\mathbf{j}i} = \prod_{l=1}^{k} \left(x_{i,j_l} - \hat{\mu}^{(0)}_{j_l}\right)
    \label{eq:moment-matrix}
\end{equation}
evaluates one condition at the fixed points: the row of ones gives the unit mass, the $n$ rows of $\mathbf{X}$ give the rule mean, and each row of $\mathbf{E}$ gives one prescribed central entry $\mathbf{j}$ of order $k \geq 2$, taken about the mean $\hat{\boldsymbol{\mu}}^{(0)}$ of the starting rule. Since each rule moment of Eq.~\eqref{eq:rule-moments} is taken about the rule mean, which depends on the weights, the moment conditions are nonlinear in $\mathbf{w}$ even with fixed points. The rows of $\mathbf{X}$ remove this nonlinearity: they keep the rule mean equal to $\hat{\boldsymbol{\mu}}^{(0)}$ in every iterate, such that the centered products in $\mathbf{E}$ are taken about a constant point and each moment condition is linear in the weights. The right-hand side keeps every prescribed moment of the starting rule unchanged, such that every iterate is as accurate as the rule it starts from. The fixed right-hand side, however, does not fix the weights: with more points than conditions, $\mathbf{A}$ has a nontrivial null space, and every weight vector $\mathbf{w}^{(0)} + \mathbf{v}$ with $\mathbf{A}\mathbf{v} = \mathbf{0}$ and entries in $[0, 1]$ matches the same moments, among which the costs $c_i^{(\ell)}$ select the solution of each iteration. Removing the points with zero weight deletes the corresponding columns of $\mathbf{A}$, such that its null space shrinks at every iteration, until the conditions determine the weights uniquely. The solutions of the linear program have inactive weights exactly equal to zero, and the corresponding points can be removed from the rule without loss of accuracy: the algorithm is summarized in Algorithm~\ref{alg:reweighting}.

\par Since the linear program of Eq.~\eqref{eq:reweighting-lp} has one equation for each condition (the $|\mathcal{J}|$ rows of $\mathbf{X}$ and $\mathbf{E}$ and the row of ones), its solutions keep at most $N_{\mathrm{C}} = |\mathcal{J}| + 1$ nonzero weights, referred to as the fixed-point lower bound. A valid rule with at most $N_{\mathrm{C}}$ points always exists \cite{bayer2006_ProofTchakaloffsTheorem}: the reweighting reaches this number with a few linear programs, but cannot go below this bound while points are fixed. While alternative weight recombination techniques exist for this same task \cite{ryu2015_ExtensionsGaussQuadrature}, the reweighted approach is adopted here because it is reused for the greedy-elimination step.

\par The fixed-point lower bound grows quickly with the dimension $n$ and the largest order $K$ of the prescribed set. The full set contains
\begin{equation}
    |\mathcal{J}_K^{\mathrm{f}}| = \binom{n+K}{K} - 1
    \label{eq:full-set-size}
\end{equation}
entries, which grow as $n^K/K!$ for large $n$, while the marginal set contains
\begin{equation}
    |\mathcal{J}_K^{\mathrm{m}}| = n + \frac{n(n+1)}{2} + n(K-2), \qquad K \geq 2
    \label{eq:marginal-set-size}
\end{equation}
entries, which grow only quadratically with $n$ for any order. For instance, for $n = 6$ and $K = 5$ the full set contains 461 entries against the 45 of the marginal set, since the mixed entries dominate the count at the higher orders.

\subsection{Greedy Point Elimination}
\label{subsec:greedy-elimination}

\par The fixed-point lower bound $N_{\mathrm{C}}$ limits rules whose points are fixed, but smaller rules can exist when the points are free to move. In Eq.~\eqref{eq:rule-design-problem}, each free point adds $n+1$ unknowns (its coordinates and its weight): matching the $N_{\mathrm{C}}$ conditions with free points therefore may require as few as
\begin{equation}
    N_{\mathrm{e}} = \left\lceil \frac{N_{\mathrm{C}}}{n+1} \right\rceil
    \label{eq:dof-estimate}
\end{equation}
points, where $\lceil \cdot \rceil$ rounds its argument up to the nearest integer. We refer to $N_{\mathrm{e}}$ as the \ac{dof} estimate: it only compares the numbers of unknowns and conditions, and it is therefore a reference value rather than a bound. A lower bound on the number of points follows instead from three requirements: $n+1$ points are needed to reproduce a full-rank covariance \cite{julier2003_SphericalSimplexUnscented}, $\lfloor K/2 \rfloor + 1$ distinct values of each variable are needed to reproduce its marginal moments up to order $K$, and, since a rule matching a full set is a cubature rule of degree $K$ (Sec.~\ref{subsec:moment-matching-cubature}), a positive rule that integrates exactly every polynomial up to degree $2s$ requires at least $\binom{n+s}{n}$ points \cite{stroud1971_ApproximateCalculationMultiple}. The combination of these requirements yields
\begin{equation}
    N_{\mathrm{lb}} =
    \begin{cases}
        \max\left\{n+1, \lfloor K/2 \rfloor + 1\right\} & \text{for marginal sets} \\[2pt]
        \max\left\{n+1, \binom{n+\lfloor K/2 \rfloor}{n}\right\} & \text{for full sets}
    \end{cases}
    \label{eq:lower-bound}
\end{equation}
where $\lfloor \cdot \rfloor$ rounds its argument down to the nearest integer (e.g., $\lfloor 5/2 \rfloor = 2$). Greedy elimination searches for rules below $N_{\mathrm{C}}$ by removing one point at a time and re-solving Eq.~\eqref{eq:rule-design-problem} for point locations and weights.

\par Each re-solve starts from the rule left by the removal of one point, with re-normalized weights and all points free. With $J = 0$, however, the solver accepts any rule that satisfies the constraints, including rules in which a point with a tiny weight lies far in the tails of the target distribution. Since the rule moment entries in Eq.~\eqref{eq:rule-moments} compute the contribution of each point only through its weight times powers of its coordinates, a very low weight can ``hide'' a distant point from every condition: in the context of uncertainty propagation and \ac{ug}, this point can introduce non-physical trajectories or inject numerical stiffness into the problem, despite its small weight. A sparsity-promoting objective in the nonlinear re-solves would favor exactly such points, as a weight can approach zero while its corresponding point drifts far into the distribution's tail. The nonlinear re-solves after the greedy elimination of a point therefore use the compactness objective
\begin{equation}
    J_{\mathrm{c}}(\mathbf{X}, \mathbf{w}) = \frac{\sum_{i=1}^{N} w_i \left\|\mathbf{z}_i\right\|^{q}}{\sum_{i=1}^{N} \bar{w}_i \left\|\bar{\mathbf{z}}_i\right\|^{q}}, \qquad q = 2\left\lceil \frac{K+1}{2} \right\rceil
    \label{eq:compactness-objective}
\end{equation}
where $\mathbf{z}_i$ are the whitened points of Eq.~\eqref{eq:whitening}, $\bar{\mathbf{z}}_i$ and $\bar{w}_i$ are those of the reduced rule from which the re-solve starts, and $q$ is the smallest even order above $K$. While the denominator serves as a normalization factor at the start of each nonlinear solve, the numerator is the rule estimate of $\mathbb{E}[\|\mathbf{z}\|^{q}]$: penalizing an even moment of order above the prescribed ones discourages the optimization from exploring tail-heavy rules with very small weights, while still satisfying the prescribed moment matching conditions. A nonlinear solve is followed by the reweighted $L_1$ minimization with fixed points of Alg.~\ref{alg:reweighting}, such that any redundant weight is removed before attempting the following greedy point elimination. 

\par
\begin{figure*}[htb]
    \centering
    \begin{subfigure}{3.25in}
        \centering
        \includegraphics[width=\linewidth]{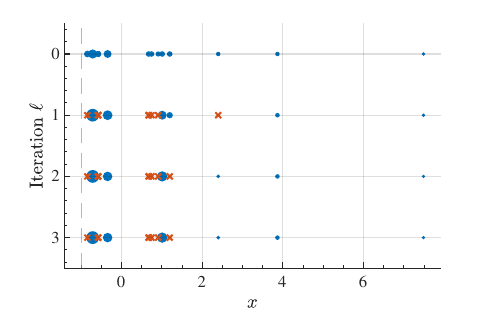}
        \caption{}
        \label{subfig:sparsification-reweighting}
    \end{subfigure}%
    \hfill
    \begin{subfigure}{3.25in}
        \centering
        \includegraphics[width=\linewidth]{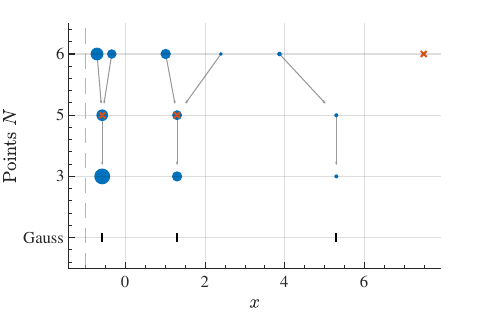}
        \caption{}
        \label{subfig:sparsification-greedy}
    \end{subfigure}
    \caption{Sparsification of a one-dimensional rule (blue marker area: weight; orange: eliminated points): a) reweighting with fixed points; b) greedy elimination with free points.}
    \label{fig:sparsification-example}
\end{figure*}

 The full sparsification algorithm is summarized in Alg.~\ref{alg:sparsification-procedure}, while Fig.~\ref{fig:sparsification-example} walks through its steps for an example one-dimensional exponential distribution, starting from $N=15$ and enforcing the moments set $\mathcal{J}_5^{\mathrm{f}}$. Figure~\ref{subfig:sparsification-reweighting} shows the iterative reweighting and weight elimination process of Sec.~\ref{subsec:reweighted-l1}: as the $\ell$ iterations advance (i.e., from the top to the bottom), some weights are pushed to 0 and eliminated (orange crosses), while points that remain active adjust their weights accordingly (indicated by the blue marker area), finishing exactly at the bound of $N_{\mathrm{C}}=|\mathcal{J}_5^{\mathrm{f}}|+1=6$ points for this one-dimensional case. Figure~\ref{subfig:sparsification-greedy} illustrates the greedy elimination procedure with free points: for the analyzed case, the first greedy removal step eliminates one term and places all points onto the three nodes of the Gaussian quadrature rule for the target, with two pairs of coincident points. The following reweighting merges each pair in a single iteration, as two coincident points give identical columns of $\mathbf{A}$ in Eq.~\eqref{eq:moment-matrix}, of which the solution of the linear program keeps at most one. We note that the final rule corresponds to the three-point Gauss rule (i.e., the positive quadrature rule with $N$ points with the highest degree of exactness of $2N-1$) computed using the Golub-Welsch algorithm, hence providing a numerical validation of the \ac{uut} solution quality for one-dimensional cases. In higher dimensions no such rule is generally available, which is why the \ac{uut} construction is carried out numerically.

\begin{algorithm}[hbt!]
\caption{Sparsification procedure}
\label{alg:sparsification-procedure}
\begin{algorithmic}[1]
\REQUIRE prescribed set $\mathcal{J}$, support $[\mathbf{a}, \mathbf{b}]$, oversized budget $N_0 > N_{\mathrm{C}}$
\ENSURE verified rule $(\mathbf{X}, \mathbf{w})$ with $N \leq N_{\mathrm{C}}$ points
\STATE Solve Eq.~\eqref{eq:rule-design-problem} with $N = N_0$ and $J = 0$, and store the verified rule as $(\mathbf{X}, \mathbf{w})$
\STATE $(\mathbf{X}, \mathbf{w}) \leftarrow \textsc{Reweight}(\mathbf{X}, \mathbf{w})$ \COMMENT{fixed points, Algorithm~\ref{alg:reweighting}}
\WHILE{$N > n+1$}
    \FOR{each of the three points with the smallest weights}
        \STATE Remove the point and rescale the remaining weights to unit mass
        \STATE Solve Eq.~\eqref{eq:rule-design-problem} with $J = J_{\mathrm{c}}$ and all points free, retrying up to twice from perturbed points \COMMENT{free points, Sec.~\ref{subsec:greedy-elimination}}
    \ENDFOR
    \IF{no candidate returns a verified rule}
        \RETURN $(\mathbf{X}, \mathbf{w})$
    \ENDIF
    \STATE $(\mathbf{X}, \mathbf{w}) \leftarrow \textsc{Reweight}(\text{most compact verified candidate})$
\ENDWHILE
\RETURN $(\mathbf{X}, \mathbf{w})$ \COMMENT{$n+1$ points reached}
\end{algorithmic}
\end{algorithm}

\begin{algorithm}[hbt!]
\caption{Reweighting with fixed points, $\textsc{Reweight}(\mathbf{X}, \mathbf{w})$}
\label{alg:reweighting}
\begin{algorithmic}[1]
\REQUIRE verified rule $(\mathbf{X}, \mathbf{w})$
\ENSURE last verified iterate $(\mathbf{X}, \mathbf{w})$, with the points of $\mathbf{X}$ unchanged
\FOR{$\ell = 0, 1, \dots$ (at most five iterations)}
    \STATE Compute the costs $c_i^{(\ell)}$ of Eq.~\eqref{eq:reweighting-costs} and solve the linear program of Eq.~\eqref{eq:reweighting-lp}
    \STATE Remove the points with zero weight and recompute the remaining weights from the conditions
    \IF{the iterate is not verified}
        \RETURN $(\mathbf{X}, \mathbf{w})$
    \ENDIF
    \STATE Store the iterate as $(\mathbf{X}, \mathbf{w})$
    \IF{the weights changed by less than the tolerance}
        \RETURN $(\mathbf{X}, \mathbf{w})$
    \ENDIF
\ENDFOR
\RETURN $(\mathbf{X}, \mathbf{w})$
\end{algorithmic}
\end{algorithm}

\section{Numerical Assessment}
\label{sec:numerical-assessment}

\par This section assesses the \ac{uut} and its sparsification against the existing unscented methods through three questions: whether the \ac{uut} recovers the existing rules where they are verified (and what it offers where they are not), how close the sparsification of Sec.~\ref{sec:rule-sparsification} comes to the smallest verified rules, and how the number of points and the matched moments scale with the number of uncertain variables. After the experimental setup, Secs.~\ref{subsec:existing-rules} to~\ref{subsec:dimension-scaling} answer the three questions in this order. 

\subsection{Experimental Setup}
\label{subsec:experimental-setup}

\par The target distributions used in this study are summarized in Table~\ref{tab:targets}, and their moments are computed analytically as detailed in Appendix~\ref{subsec:target-moments}, such that no approximations can affect the comparison. All the performed experiments include four ``core distributions'': Gaussian (on which most of the reference unscented methods are designed), uniform (bounded and symmetric), beta (bounded and skewed), and the exponential distribution (bounded below and noticeably skewed). Multivariate targets have independent components, except for the correlated Gaussian targets (used in Sec.~\ref{subsec:point-efficiency}) and the mixed exponentials of Sec.~\ref{subsec:dimension-scaling}. The lognormal, Weibull, and Gaussian mixture families, together with the parameter sweeps of Table~\ref{tab:targets}, enter only the aggregated results of Secs.~\ref{subsec:existing-rules} and~\ref{subsec:point-efficiency}. While the core families cover bounded and unbounded, symmetric and skewed targets, these extend the comparison to heavy tails, negative skewness, and bimodal targets. The prescribed sets are the full and marginal sets of Sec.~\ref{subsec:moment-matching-cubature} and their unions.

\begin{table}[hbt!]
\caption{\label{tab:targets} Target distributions (bold: parameter of the core targets; the gamma law with $\alpha = 1$ is the exponential). All distributions are standardized (i.e., zero mean and unit variance).}
\centering
\setlength{\tabcolsep}{4.5pt}
\begin{tabular}{llllcc}
\hline\hline
Family & Law & Parameter values & Support & Skewness & Kurtosis \\
\hline
Gaussian & $\mathcal{N}(0, 1)$ & -- & $\mathbb{R}$ & 0 & 3 \\
Uniform & $\mathcal{U}(0, 1)$ & -- & $[0, 1]$ & 0 & 1.8 \\
Beta & $\mathrm{Beta}(2, b)$ & $b$: 2, 3, 6, \textbf{12}, 24, 48 & $[0, 1]$ & 0 to 1.29 & 2.14 to 5.33 \\
Gamma & $\mathrm{Gamma}(\alpha, 1)$ & $\alpha$: 0.5, \textbf{1}, 2, 4, 8, 16 & $[0, \infty)$ & 0.50 to 2.83 & 3.38 to 15.0 \\
\hline
Correlated Gaussian\textsuperscript{a} & $\mathcal{N}(\mathbf{0}, \mathbf{R}_{\rho})$ & $\rho$: 0.5, 0.9, 0.99 & $\mathbb{R}^n$ & 0 & 3 \\
Mixed exponentials\textsuperscript{b} & $\mathbf{R}_{\rho}^{1/2}\boldsymbol{\epsilon}$ & $\rho$: 0.5 & Shifted cone & 1.36 to 1.84 & 6.33 to 8.25 \\
\hline
Lognormal & $\mathrm{LogN}(0, s^2)$ & $s$: 0.1, 0.25, 0.5, 0.7, 0.75, 1 & $(0, \infty)$ & 0.30 to 6.19 & 3.16 to 114 \\
Weibull & $\mathrm{Weibull}(c, 1)$ & $c$: 0.75, 1, 1.5, 2, 3.6, 5 & $[0, \infty)$ & $-0.25$ to 3.12 & 2.72 to 19.0 \\
Gaussian mixture & $\frac{1}{2}\mathcal{N}(\pm d, \frac{1}{4})$ & $d$: 0.5, 1, 1.5, 2, 3, 4 & $\mathbb{R}$ & 0 & 1.06 to 2.50 \\
\hline\hline
\multicolumn{6}{p{\dimexpr\linewidth-2\tabcolsep\relax}}{\textsuperscript{a} $\mathbf{R}_{\rho} = (1 - \rho)\mathbf{I} + \rho\mathbf{1}\mathbf{1}^{T}$ is the equicorrelation matrix, and $\mathbf{R}_{\rho}^{1/2}$ its matrix square root. The targets of the first and last groups have independent components.} \\
\multicolumn{6}{p{\dimexpr\linewidth-2\tabcolsep\relax}}{\textsuperscript{b} $\boldsymbol{\epsilon}$ has $n$ independent standardized exponential components, such that $\mathbf{R}_{\rho}^{1/2}\boldsymbol{\epsilon}$ has covariance $\mathbf{R}_{\rho}$ and nonzero mixed third-order entries (e.g., 0.61 for $n = 2$). No bound is enforced on its support, to keep the support constraints of Eq.~\eqref{eq:valid-rule} as a box. Skewness and kurtosis are marginal values, from $n = 12$ to $n = 2$.}
\end{tabular}
\end{table}

\par We choose a small set of reference rules that is representative of existing techniques to construct sigma-point rules: the \ac{ut} \cite{julier2000_NewMethodNonlinear}, the \ac{gut} \cite{ebeigbe2021_GeneralizedUnscentedTransformation}, and the \ac{cut4} \cite{adurthi2012_ConjugateUnscentedTransform, adurthi2013_ConjugateUnscentedTransform}. The \ac{ut} uses the scale $\sqrt{n+1}$, meaning that its central weight is always $1/(n+1)$ and stays positive, such that the \ac{ut} rule can violate the validity conditions only through points outside the support \cite{ebeigbe2021_GeneralizedUnscentedTransformation}. The \ac{gut} receives the support bounds of the target distribution. The \ac{cut4} uses its uniform construction for the non-scalar uniform targets and its Gaussian construction for every other distribution: on these targets, the \ac{cut4} is a transferred baseline, since no \ac{cut4} exists for them. Each reference goes through the same verification as the \ac{uut} rules on the prescribed set they are built to match, as illustrated in Table~\ref{tab:existing-rules}: $\mathcal{J}_4^{\mathrm{m}}$ for the \ac{gut} on any target, and $\mathcal{J}_3^{\mathrm{f}}$ and $\mathcal{J}_5^{\mathrm{f}}$ for the \ac{ut} and the \ac{cut4}, respectively. All \ac{uut} solves use the moment band $\varepsilon_m = 10^{-8}$ of Eq.~\eqref{eq:rule-design-problem} and \ac{ipopt} tolerance $10^{-8}$, the reweighting stops when the weights change by less than $10^{-6}$ between iterations, and every rule (including the ones built from reference methods) is independently verified with $\varepsilon_v = 10^{-6}$. The \ac{uut} source code is available at \url{https://github.com/RiccardoMinnozzi-UniBo/Unified-Unscented-Transform}

\subsection{Existing Rules as Particular Solutions}
\label{subsec:existing-rules}

\par
\begin{figure*}[htb]
    \centering
    \begin{subfigure}{3.25in}
        \centering
        \includegraphics[width=\linewidth]{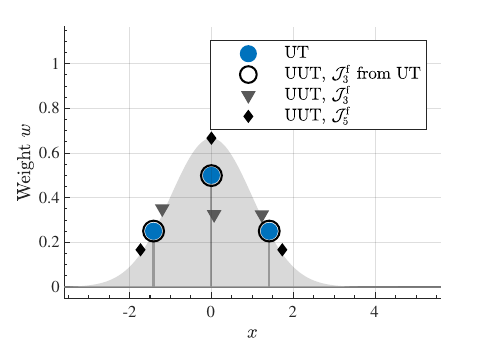}
        \caption{}
        \label{subfig:rules-ut}
    \end{subfigure}%
    \hfill
    \begin{subfigure}{3.25in}
        \centering
        \includegraphics[width=\linewidth]{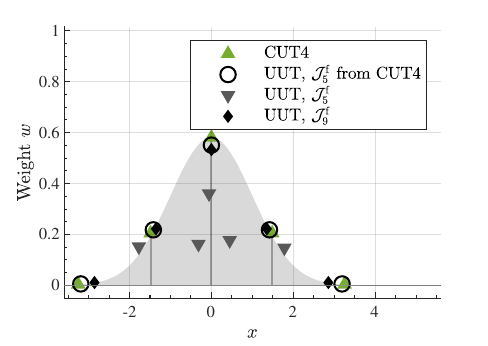}
        \caption{}
        \label{subfig:rules-cut4}
    \end{subfigure}

    \vspace{0.5\baselineskip}
    \begin{subfigure}{3.25in}
        \centering
        \includegraphics[width=\linewidth]{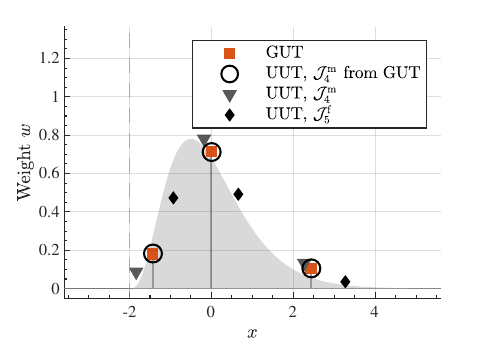}
        \caption{}
        \label{subfig:rules-gut}
    \end{subfigure}%
    \hfill
    \begin{subfigure}{3.25in}
        \centering
        \includegraphics[width=\linewidth]{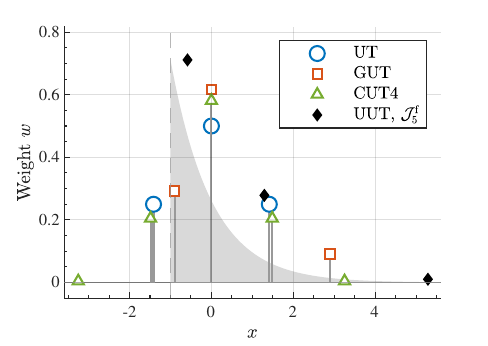}
        \caption{}
        \label{subfig:rules-exponential}
    \end{subfigure}
    \caption{One-dimensional rules, with each point at its weight over the target density: a) \ac{ut}, Gaussian; b) \ac{cut4}, Gaussian; c) \ac{gut}, gamma; d) exponential (the empty markers indicate failed verification of the rule).}
    \label{fig:existing-rules}
\end{figure*}

 Figure~\ref{fig:existing-rules} compares the rules in one dimension, allowing intuitive visualization of the equivalence between the \ac{uut} and existing methods: while similar conclusions can be derived from higher-dimensional test cases, the one-dimensional comparison is preferred for its intuitiveness and the availability of analytically verified quadrature rules for independent verification. Figures~\ref{subfig:rules-ut},~\ref{subfig:rules-cut4}, and~\ref{subfig:rules-gut} show each reference on a target where it is verified, together with three \ac{uut} rules built with the same $N$: the rule matching the prescribed set of the reference initialized from the reference itself (the ring around each reference point), the one obtained from the default \ac{uut} initial guess, and the unique rule matching $\mathcal{J}_{2N-1}^{\mathrm{f}}$ (i.e., the \ac{uut} construction of the $N$-point Gauss rule of the target). The \ac{ut} (Fig.~\ref{subfig:rules-ut}) and \ac{cut4} (Fig.~\ref{subfig:rules-cut4}) are shown on the Gaussian target they are designed for, while the \ac{gut} (Fig.~\ref{subfig:rules-gut}) is shown on the gamma target with $\alpha = 4$, a skewed target on which it is verified.

\par For the targets in Figs.~\ref{subfig:rules-ut},~\ref{subfig:rules-cut4}, and~\ref{subfig:rules-gut}, the \ac{uut} rules initialized from the references are observed to closely or exactly match the provided guess, showing that existing sigma-point schemes can be interpreted as particular solutions to the \ac{uut} construction. However, when using the default \ac{uut} initial guess, different rules are identified in all three cases: since all the accepted rules are independently checked (i.e., they match the prescribed moments up to tolerance $\varepsilon_v$), this provides a numerical verification that the references leave free \acp{dof} in one dimension. None of their prescribed moment sets, in fact, reaches the degree of exactness $2N-1$ of the Gauss rule, which instead guarantees the uniqueness of the identified rule. The same features are observed when repeating the analysis with 2 to 8 variables: the \ac{uut} closely matches each verified reference when started from it, while it returns different verified rules when started from the default initial guess. The \ac{uut} rules built from the prescribed $\mathcal{J}_{2N-1}^{\mathrm{f}}$ sets match the corresponding Gauss rules for all one-dimensional targets, including those where all the existing references fail, such as the exponential distribution shown in Fig.~\ref{subfig:rules-exponential}.

\par While Fig.~\ref{fig:existing-rules} provides visual insight on the performance and failure modes of existing rules on a few representative targets, performing the analysis on a wider range of distributions (and of dimensions, as done in Sec.~\ref{subsec:dimension-scaling}) exposes additional failure mechanisms. These failures are of two kinds: either a reference is applied to a target it is not designed for (e.g., the \ac{cut4} transferred to skewed targets), or its structure violates the validity conditions of Eq.~\eqref{eq:valid-rule} on the targets it is designed for (e.g., the negative central weight of the \ac{gut} for Gaussian targets with more than three variables). The \ac{ut} places a point outside of the support of several bounded skewed targets (two each among the gamma, lognormal, and Weibull families), while the \ac{gut} satisfies the support-set constraint at the cost of part of its moment-matching capability. With more variables, the \ac{gut} central weight becomes negative for the independent Gaussian targets with more than three variables and for all the non-scalar uniform targets, while the \ac{ut} points, placed at $\sqrt{n+1}\sigma$ from the mean to prevent negative weights, tend to leave tight supports as the dimension grows (e.g., for the uniform targets from three variables, where $\sqrt{n+1} > \sqrt{3}$). The \ac{cut4} is verified on the Gaussian and uniform targets it is designed for, while transferring its Gaussian construction to other distributions fails because of their tight supports, nonzero odd moments, or fourth-order moments that differ from the Gaussian ones. The \ac{uut}, instead, returns a verified rule for $\mathcal{J}_5^{\mathrm{f}}$ on every scalar target (verified against Golub-Welsch quadrature construction), and every \ac{uut} rule satisfies the validity conditions of Eq.~\eqref{eq:valid-rule} by construction.

\subsection{Point Efficiency and Sparsification}
\label{subsec:point-efficiency}

\par We now evaluate the effectiveness of the sparsification technique presented in Sec.~\ref{sec:rule-sparsification} against a simple direct search, which solves Eq.~\eqref{eq:rule-design-problem} from the default \ac{uut} initial guess at a sequence of fixed budgets and reports the smallest number of points $N$ that yields a verified rule. The two approaches explore the point budget from opposite directions: the sparsification works top-down from an oversized rule, while the direct search works bottom-up from the estimate $N_{\mathrm{e}}$: the  agreement of the two approaches provides an indicator for how reliable and effective they are in approximating the smallest sigma-point rule.

\par Both approaches are assessed on two-dimensional distributions with full sets of orders three to six: the Gaussian, uniform, beta, exponential, Gaussian mixture, and lognormal targets, and the correlated Gaussian targets with correlation 0.5 and 0.9. For each case, the direct search tries the budgets from $N_{\mathrm{e}}-2$ to $N_{\mathrm{e}}+3$, never below $N_{\mathrm{lb}}$, stopping at the first one that yields a verified rule. The sparsification, instead, starts from an oversized rule of $\lceil 1.5\,N_{\mathrm{C}} \rceil$ points and follows Algorithm~\ref{alg:sparsification-procedure}.

\begin{figure*}[htb]
    \centering
    \begin{subfigure}{3.25in}
        \centering
        \includegraphics[width=\linewidth]{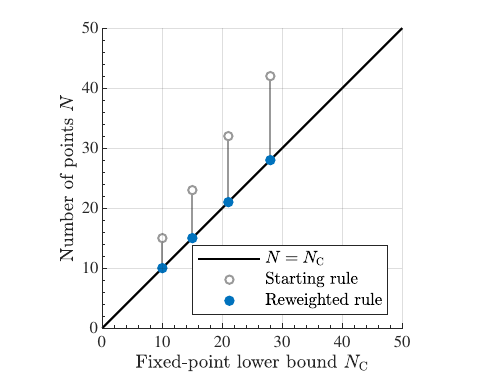}
        \caption{}
        \label{subfig:reweighting-bound}
    \end{subfigure}%
    \hfill
    \begin{subfigure}{3.25in}
        \centering
        \includegraphics[width=\linewidth]{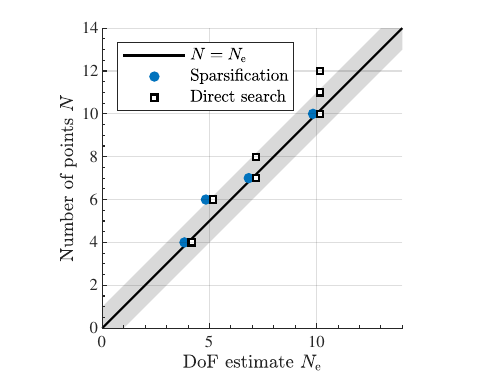}
        \caption{}
        \label{subfig:sparsification-estimate}
    \end{subfigure}
    \caption{Sparsification of the benchmark cases: a) reweighting from $\lceil 1.5\,N_{\mathrm{C}} \rceil$ points against $N_{\mathrm{C}}$; b) sparsified rules and direct search against $N_{\mathrm{e}}$.}
    \label{fig:point-efficiency}
\end{figure*}

\par Figure~\ref{fig:point-efficiency} separates the contributions of the two stages of the sparsification. Figure~\ref{subfig:reweighting-bound} compares the rules left by the fixed-point reweighting of Sec.~\ref{subsec:reweighted-l1} with its theoretical lower bound $N_{\mathrm{C}}$, and Fig.~\ref{subfig:sparsification-estimate} compares the rules of the complete sparsification and of the direct search with the \ac{dof} estimate of Eq.~\eqref{eq:dof-estimate}. In both figures the blue markers represent the result for each distribution and prescribed set: the overlap of all blue markers in both figures indicates that both steps of the sparsification process identify rules with the same sparsity across the analyzed distributions and moment sets. More noticeably, Fig.~\ref{subfig:sparsification-estimate} shows that both the direct search and the greedy point elimination identify rules close to the \ac{dof} estimate $N_e$, and that the greedy elimination always finds solutions that are at least as sparse as what is found with the direct search, with the added benefit that the sparsification procedure promotes compact rules, avoiding the small weights associated to points far away from the mean that are instead observed in direct solves.

\subsection{Dimension Scaling}
\label{subsec:dimension-scaling}

\par We now assess how the \ac{uut} scales with the number of uncertain variables through two (closely related) questions: whether it matches the moments of each reference with fewer points, and whether it matches more moments with the same points. Both comparisons take the point count of each reference as the budget of the \ac{uut}, such that the \ac{ut} and \ac{gut} comparisons span 2 to 12 variables ($N = 2n+1$), while the \ac{cut4} comparisons stop at 6 variables, since the Gaussian construction that sets their point count grows as $2n + 2^n$ (which would amount to $N=142$ for $n=7$). For the first question, the prescribed set of each reference is initialized at its point budget and reduced through Algorithm~\ref{alg:sparsification-procedure}. For the second question, Eq.~\eqref{eq:rule-design-problem} is solved at the reference point count for moment sets that grow by one marginal order at a time on top of the full set (i.e., $\mathcal{J}_K^{\mathrm{f}}$, then $\mathcal{J}_K^{\mathrm{f}} \cup \mathcal{J}_{K+1}^{\mathrm{m}}$, and so on up to order six), stopping at the first failure. The starting full sets range from $\mathcal{J}_2^{\mathrm{f}}$ to $\mathcal{J}_3^{\mathrm{f}}$ with the $2n+1$ points shared by the \ac{ut} and \ac{gut}, and from $\mathcal{J}_3^{\mathrm{f}}$ to $\mathcal{J}_5^{\mathrm{f}}$ with the \ac{cut4} points, keeping the largest set matched on each target. All the targets of Table~\ref{tab:targets} enter the comparison, and the figure reports the least favorable one: the largest sparsified rule among the targets where one is verified, and the largest set matched on every target.

\begin{figure*}[htb]
    \centering
    \begin{subfigure}{3.25in}
        \centering
        \includegraphics[width=\linewidth]{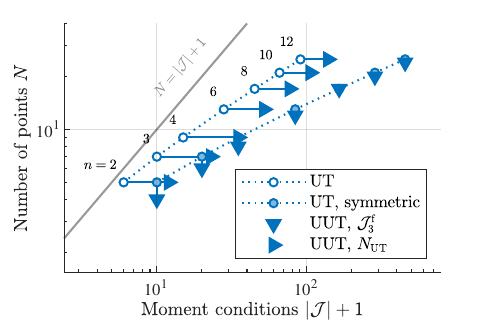}
        \caption{}
        \label{subfig:budget-ut}
    \end{subfigure}%
    \hfill
    \begin{subfigure}{3.25in}
        \centering
        \includegraphics[width=\linewidth]{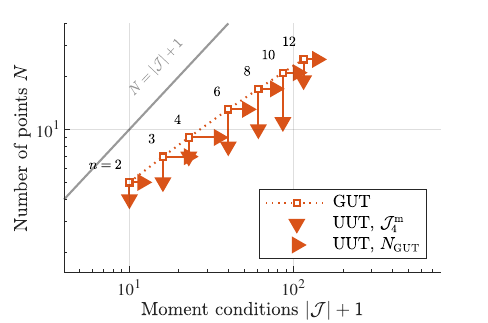}
        \caption{}
        \label{subfig:budget-gut}
    \end{subfigure}

    \vspace{0.5\baselineskip}
    \begin{subfigure}{3.25in}
        \centering
        \includegraphics[width=\linewidth]{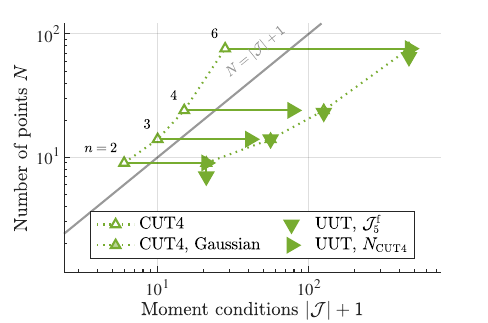}
        \caption{}
        \label{subfig:budget-cut}
    \end{subfigure}%
    \hfill
    \begin{subfigure}{3.25in}
        \centering
        \includegraphics[width=\linewidth]{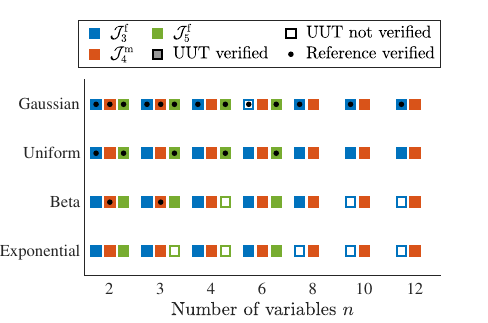}
        \caption{}
        \label{subfig:verification-map}
    \end{subfigure}
    \caption{Fewer points and more moments with 2 to 12 variables: a--c) points against moment conditions for the \ac{ut}, \ac{gut}, and \ac{cut4}; d) targets with a verified rule.}
    \label{fig:dimension-budget}
\end{figure*}

\par Figure~\ref{fig:dimension-budget} places each rule by its numbers of moment conditions and points, such that rules below the line $N = |\mathcal{J}| + 1$ have fewer points than conditions. Each reference is drawn at its prescribed set (shaded marker) and at the set it matches on a generic target (open marker): the two coincide for the \ac{gut}, while the \ac{ut} and \ac{cut4} match only $\mathcal{J}_2^{\mathrm{f}}$ outside the symmetric and Gaussian targets, respectively. These markers report the theoretical features of each reference (i.e., its number of points and prescribed set when its assumptions are verified) regardless of the failures discussed in Sec.~\ref{subsec:existing-rules}, while every \ac{uut} rule is valid by construction. The downward-pointing triangles are the sparsified \ac{uut} rules matching the prescribed set of the reference (i.e., vertical distance indicates the point reduction), while the right-pointing triangles are the largest sets matched with the same points (i.e., horizontal distance from the open marker indicates the additional moment conditions). Figure~\ref{subfig:verification-map} shows on which targets each rule is verified on the prescribed set of the reference (squares for the sparsified rules, dots for the references), where an empty square means that the sparsification found no verified rule.

\par Figures~\ref{subfig:budget-ut}--\ref{subfig:budget-cut} show that the \ac{uut} matches or outperforms the existing rules on the whole range of targets, while ensuring that each rule is a valid probability measure. The comparison is conservative for the \ac{uut}, as the reference markers ignore the failures of the existing rules discussed in Sec.~\ref{subsec:existing-rules} and shown in Fig.~\ref{subfig:verification-map}. Since the \ac{gut} is designed to match its set on any target, its shaded and open markers coincide: Fig.~\ref{subfig:budget-gut} provides the most straightforward comparison, as the \ac{uut} rules both use fewer points for the same moment conditions (i.e., downward-pointing triangles below the \ac{gut} markers) and match more moment conditions with the same points (i.e., right-pointing triangles to the right of them). The \ac{ut} and \ac{cut4} comparisons lead to the same conclusions on generic targets, where the \ac{uut} matches their prescribed sets with fewer points and more than $\mathcal{J}_2^{\mathrm{f}}$ with the same points. On the targets these references are designed for, the least favorable \ac{uut} rules fall short of their prescribed sets as the variables grow: these rules correspond to skewed targets, however, and an a-posteriori check reveals that the \ac{uut} still matches at least the same moments with fewer points.

\par Overall, the comparisons of Fig.~\ref{fig:dimension-budget} address the two limitations of the existing rules noted in Sec.~\ref{sec:introduction}: fewer points for the same moments give a smaller guidance problem with the same statistics, while more moments for the same points give more accurate statistics at the same problem size. Unlike the references, which are verified only on the few targets they are designed for, the \ac{uut} offers both sides on most of the distributions tested here, without assumptions on the target beyond its moments and support.

\section{Application to Unscented Guidance}
\label{sec:unscented-guidance}

\par Sections~\ref{sec:unified-unscented-transform} to~\ref{sec:numerical-assessment} assessed the \ac{uut} rules as cubature schemes, through their validity, number of points, and matched moments. This section assesses the \ac{uut} capability to integrate within the \ac{ug} framework of Refs.~\cite{ross2015_RiemannStieltjesOptimal, ross2025_UnscentedTrajectoryOptimization}. The section first presents an overview of Unscented Optimal Control and the role that cubature rules play inside it, then provides some insights on the enforcement of chance constraints with the available statistical information, and finally showcases the overall framework applied to an uncertain atmospheric entry problem.

\subsection{Unscented Optimal Control}
\label{subsec:unscented-optimal-control}

\par Unscented Optimal Control addresses tychastic (i.e., conditionally deterministic) systems, meaning systems subject to parametric uncertainties for which, if the parameter values are known exactly, we have a fully-deterministic model available \cite{ross2014_UnscentedOptimalControl, ross2016_PathConstraintsTychastic}: the technique aims to identify optimal controls over the many possible realizations of the uncertain parameters. Objectives and constraints therefore represent expectations or quantiles over the parameter distributions: these expectations are formalized as Riemann-Stieltjes integrals \cite{ross2015_RiemannStieltjesOptimal} and approximated using cubature rules over the distribution support. While any cubature rule is a valid approach, unscented methods provide a straightforward technique to build cheap approximations of the expectation integrals through their sigma-points and associated weights (as long as the cubature rule provides a valid probability measure as discussed in previous sections). The \ac{ocp} statistics are propagated using the unscented approach, defining a copy of the state trajectory for each sigma point and reconstructing the desired quantities through the weighted sum of Eq.~\eqref{eq:sigma-point-rule}.

\par For conditionally-deterministic \acp{ocp} with arbitrary costs and constraints, a single control policy may not always be able to satisfy the original problem formulation, thus multiple reformulations of the same unscented \ac{ocp} can be adopted to target different aspects of the original \ac{ocp} (the interested reader is referred to Refs.~\cite{ross2014_UnscentedOptimalControl, ross2016_PathConstraintsTychastic, ross2025_UnscentedTrajectoryOptimization} for further details). For the scope of this manuscript, we present one potential instance of unscented \ac{ocp} that matches that of the example application considered later. Consider a system with state $\mathbf{y}(t) \in \mathbb{R}^{n_y}$ and control $\mathbf{u}(t) \in \mathbb{R}^{n_u}$, whose dynamics $\dot{\mathbf{y}} = \mathbf{f}(\mathbf{y}, \mathbf{u}, \mathbf{x})$ depend on the uncertain variables $\mathbf{x}$ of Sec.~\ref{subsec:moment-matching-cubature}. Given a rule $(\mathbf{X}, \mathbf{w})$ for their distribution, each sigma point $\mathbf{x}_i$ defines a copy of the trajectory $\mathbf{y}_i(t)$, obtained by integrating the dynamics from the initial state $\mathbf{y}_i(t_0)$ under the common control input $\mathbf{u}(t)$ (note that the common control input need not be a feedforward control in the classic sense \cite{han2025_UnscentedGuidanceFeedback}). For an \ac{ocp} with Mayer cost $\varphi(\mathbf{y}_i(t_f), t_f, \mathbf{x}_i)$, endpoint constraints $\mathbf{e}(\mathbf{y}_i(t_f), t_f, \mathbf{x}_i)=\mathbf{e}_f$ and path constraints $c(\mathbf{y}_i, \mathbf{u}, t, \mathbf{x}_i)\leq c_{\max}$, the unscented \ac{ocp} may be formulated as
\begin{equation}
    \begin{aligned}
        \min_{\mathbf{u},\, t_f} \quad & \sum_{i=1}^{N} w_i\, \varphi\left(\mathbf{y}_i(t_f), t_f, \mathbf{x}_i\right) \\
        \text{s.t.} \quad & \dot{\mathbf{y}}_i = \mathbf{f}\left(\mathbf{y}_i, \mathbf{u}, \mathbf{x}_i\right), & & i = 1, \dots, N \\
        & \sum_{i=1}^{N} w_i\, \mathbf{e}\left(\mathbf{y}_i(t_f), t_f, \mathbf{x}_i\right) = \mathbf{e}_f \\
        & \Pr[c(\mathbf{y}, \mathbf{u}, t, \mathbf{x})\leq c_{\max}] \geq 1-\varepsilon, & & t \in [t_0, t_f]
    \end{aligned}
    \label{eq:unscented-ocp}
\end{equation}
where the last row relating to the path constraints is expressed in chance-constraint form, for which a deterministic reformulation is discussed in more detail in Sec.~\ref{subsec:moment-based-risk-constraint}. The particular unscented \ac{ocp} formulation provided in Eq.~\eqref{eq:unscented-ocp} is one potential approach to formulate a tychastic optimal control problem: the objective represents the mean Mayer cost, the endpoint constraints are also enforced on the mean, while the path constraints are enforced up to an acceptable risk fraction $\varepsilon$. While this formulation is shown here as it matches the example problem in Sec.~\ref{subsec:entry-problem}, it is important to understand that the cubature approximations used for the costs, endpoint, and path constraints can be applied to any quantity in a generic \ac{ocp} formulation, examples of which include problems with uncertain initial conditions (where $\mathbf{x}_i=\mathbf{y}_i(t_0)$) and covariance steering problems (where the uncertain initial conditions are combined with an endpoint constraint on the covariance of the states, $\sum_{i=1}^{N} w_i \left(\mathbf{y}_i(t_f) - \bar{\mathbf{y}}_f\right)\left(\mathbf{y}_i(t_f) - \bar{\mathbf{y}}_f\right)^T \preceq \mathbf{P}_f$, with mean $\bar{\mathbf{y}}_f = \sum_{i=1}^{N} w_i \mathbf{y}_i(t_f)$ and target covariance $\mathbf{P}_f$). The unscented cubature rules provide a generic deterministic reformulation of the uncertain quantities in the \ac{ocp}, making them compatible with arbitrary transcription techniques, such as the pseudospectral methods implemented in the \acs{raptor} solver used in this study \cite{minnozzi2026_TheoreticalComputationalAdvances}.

\subsection{Moment-Based Chance Constraints}
\label{subsec:moment-based-risk-constraint}

\par While the approximation of expectation integrals through unscented methods follows from applying Eq.~\eqref{eq:sigma-point-rule}, handling chance constraints requires a way to approximate the probability of satisfying the specific constraint function under uncertainty. Since sigma-point rules directly approximate moments but not probabilities, alternative methods are required to handle this class of constraints. A workaround is proposed in Ref.~\cite{schwartz2025_UnscentedHypersonicTrajectory}, adopting a "deviation around the mean" approach: given the constraint $c$, Eq.~\eqref{eq:sigma-point-rule} is used to compute its mean $\mu_c$ and standard deviation $\sigma_c$ (i.e., the square root of the variance $\mathbb{E}[(c-\mu_c)^2]$), and the chance constraint is enforced as
\begin{equation}
    \mu_c+\kappa\sigma_c\leq c_{\max}
    \label{eq:kappa-chance-constraint}
\end{equation}
where the parameter $\kappa$ controls the desired confidence in satisfying the constraint $c\leq c_{\max}$. Ref.~\cite{schwartz2025_UnscentedHypersonicTrajectory} interprets the value $\kappa$ as the number of standard deviations of confidence for the satisfaction of the chance constraint: choosing $\kappa=3$ means targeting a 3-sigma confidence (i.e., $99.865\%$ probability of constraint satisfaction). This value is only exact if the statistics reconstructed via the unscented rule are accurate, and if the realization of the chance constraint follows a Gaussian distribution.

\par As discussed in Sec.~\ref{subsec:moment-matching-cubature}, an unscented rule that matches the prescribed set $\mathcal{J}_K^{\mathrm{f}}$ (i.e., all moments up to order $K$) is an exact cubature rule for all polynomials in the uncertain variable $\mathbf{x}$ of degree up to $K$: if $c$ depends on $\mathbf{x}$ as a polynomial of degree up to $K$, the $\mu_c$ computed through the unscented rule is exact. On the other hand, since the standard deviation $\sigma_c$ requires the expectation of $c^2$, its value will be exact only for polynomials up to degree $\lfloor K/2 \rfloor$. While these exactness considerations may not hold for arbitrary \acp{ocp}, they provide useful insights on the relevant error sources when using such chance constraint reformulations: through the Taylor expansion of $c$ about the mean of $\mathbf{x}$, the order $K$ sets which terms of the nonlinear response enter the statistics exactly. With $K = 2$ (i.e., mean and covariance only), the second-order terms enter $\mu_c$, but only the linear ones enter $\sigma_c$. The third-order entries add third-order terms to $\mu_c$, and the products of the linear and second-order terms to $\sigma_c$.

\par For constraint realizations that are distributed according to a Gaussian distribution, the constraint formulated according to Eq.~\eqref{eq:kappa-chance-constraint} implies $\Pr[c\leq c_{\max}]\geq \Phi(\kappa)$, where $\Phi(\cdot)$ is the cumulative distribution function of the standard Gaussian distribution (hence the $\Phi(3)=99.865\%$ mentioned previously). For a general \ac{ocp} under arbitrary uncertainties, however, assuming that a value $c$ follows a Gaussian distribution may not be accurate. Enforcing chance constraints in a meaningful way requires being able to precisely match the prescribed violation probabilities (i.e., exploiting the available design space up to the acceptable risk). More accurate representations can be obtained from higher-order statistical information: the Cornish-Fisher expansion approximates the quantile of order $1-\varepsilon$ of a distribution through its moments \cite{fisher1960_PercentilePointsDistributions, burnett2026_AnalyticNonGaussianConfidence}. Retaining the first correction only (i.e., the contribution of the skewness), the chance constraint can be reformulated as
\begin{equation}
    \mu_c(t) + z\, \sigma_c(t) + \frac{z^2 - 1}{6} \max\left(0, \frac{\mu_{3,c}(t)}{\sigma_c^2(t)}\right) \leq c_{\max}, \qquad z = \Phi^{-1}(1 - \varepsilon)
    \label{eq:cf-constraint}
\end{equation}
where $\mu_{3,c}(t) = \sum_{i=1}^{N} w_i \left(c_i(t) - \mu_c(t)\right)^3$ is the third central moment of $c$ reconstructed by the rule. Equation~\eqref{eq:cf-constraint} is the constraint of Eq.~\eqref{eq:kappa-chance-constraint} with a margin that follows the skewness predicted by the rule, and it reduces to it with $\kappa = z$ when the predicted output is symmetric. The design parameter thus becomes the allowed exceedance probability $\varepsilon$ instead of $\kappa$. The skewness term is bounded below by zero: as the expansion is valid for distributions that are close to a Gaussian, this bound prevents the optimizer from exploiting constraint distributions whose skewness may produce incorrect approximations of the associated quantiles (i.e., for strongly negative skewness, Eq.~\eqref{eq:cf-constraint} may predict a risk that decreases as the target quantile increases). Since the maximum is not differentiable at zero, the transcription replaces it with the smooth approximation $(s + \sqrt{s^2 + \delta^2})/2$ of $\max(0, s)$, with $\delta=10^{-5}$ small with respect to $\sigma_c$. 

\par For the specific case of path constraints (e.g., in the problem formulation of Eq.~\eqref{eq:unscented-ocp}), the constraints of Eqs.~\eqref{eq:kappa-chance-constraint},~\eqref{eq:cf-constraint} are enforced at each instant separately, meaning that they bound the pointwise exceedance probability $\Pr[c(t) > c_{\max}]$. The quantity of interest for the mission is instead the trajectory risk, i.e., the probability that $c$ exceeds $c_{\max}$ at least once along the trajectory, or equivalently that the peak of $c(t)$ exceeds $c_{\max}$. Since any realization that exceeds the limit at a given instant also exceeds it along the trajectory, the trajectory risk is always greater (or equal to, if all trajectories reach their peak at the same time) than the largest pointwise probability. 

\subsection{Entry Problem}
\label{subsec:entry-problem}

\par We now introduce a representative \ac{ug} problem, taken from \cite{schwartz2025_UnscentedHypersonicTrajectory}: the crossrange maximization of a reusable entry vehicle under uncertain atmospheric density. The vehicle is modeled as a point mass over a spherical, rotating Earth, with state $\mathbf{y} = [h, \theta, \phi, v, \gamma, \psi]^T$, namely the altitude, longitude, latitude, velocity, flight-path angle, and heading angle, and controls $\mathbf{u} = [\sigma, \alpha]^T$ being the bank angle $\sigma$ and the angle of attack $\alpha$. The equations of motion are
\begin{equation}
    \begin{aligned}
        \dot{h} &= v \sin\gamma \qquad
        \dot{\theta} = \frac{v \cos\gamma \cos\psi}{r \cos\phi} \qquad
        \dot{\phi} = \frac{v \cos\gamma \sin\psi}{r} \\
        \dot{v} &= -\frac{D}{m} - g \sin\gamma + \omega_{\oplus}^2 r \cos\phi \left(\sin\gamma \cos\phi - \cos\gamma \sin\phi \sin\psi\right) \\
        \dot{\gamma} &= \frac{L \cos\sigma}{m v} - \left(\frac{g}{v} - \frac{v}{r}\right) \cos\gamma + 2\omega_{\oplus} \cos\phi \cos\psi + \frac{\omega_{\oplus}^2 r}{v} \cos\phi \left(\cos\gamma \cos\phi + \sin\gamma \sin\phi \sin\psi\right) \\
        \dot{\psi} &= \frac{L \sin\sigma}{m v \cos\gamma} - \frac{v}{r} \cos\gamma \cos\psi \tan\phi + 2\omega_{\oplus} \left(\tan\gamma \cos\phi \sin\psi - \sin\phi\right) - \frac{\omega_{\oplus}^2 r}{v \cos\gamma} \sin\phi \cos\phi \cos\psi
    \end{aligned}
    \label{eq:entry-dynamics}
\end{equation}
where $r = R_{\oplus} + h$ is the distance from the center of the Earth, $R_{\oplus}$ and $\omega_{\oplus}$ are the radius and rotation rate of the Earth, $g = \mu_{\oplus}/r^2$ is the gravitational acceleration, with gravitational parameter $\mu_{\oplus}$, $m$ is the mass of the vehicle, and $L$ and $D$ are the lift and drag forces.

\par The lift and drag forces are computed as $L = \frac{1}{2}\rho v^2 S C_L$ and $D = \frac{1}{2}\rho v^2 S C_D$, where $S$ is the reference area and the coefficients $C_L = c_{L,0} + c_{L,1}\alpha$ and $C_D = c_{D,0} + c_{D,1}\alpha + c_{D,2}\alpha^2$ depend on the angle of attack only. The density follows the exponential model
\begin{equation}
    \rho(h) = \rho_0 e^{-\beta h}
    \label{eq:density-model}
\end{equation}
where $\rho_0$ is the density at zero altitude and $\beta$ is the inverse scale height. The heating rate is
\begin{equation}
    q = c_q \left(c_0 + c_1\alpha + c_2\alpha^2 + c_3\alpha^3\right) \rho^{1/2}\, v^{3.07}
    \label{eq:heating-rate}
\end{equation}
where $c_q$ and $c_0, \dots, c_3$ are the coefficients of the heating model. Table~\ref{tab:entry-data} summarizes all relevant quantities.

\par The vehicle starts from the conditions of Table~\ref{tab:entry-data}, and the objective is the final latitude $\phi(t_f)$: since the vehicle turns towards negative latitudes, minimizing $\phi(t_f)$ maximizes the crossrange, as in Ref.~\cite{schwartz2025_UnscentedHypersonicTrajectory}. The problem formulation follows that of the unscented \ac{ocp} in Eq.~\eqref{eq:unscented-ocp}: the mean terminal crossrange is maximized (i.e., $\varphi=\phi(t_f)$), endpoint constraints are enforced only on their mean (i.e., $\mathbf{e} = [h, v, \gamma]^T$), and the path constraint enforced in probabilistic terms is the heat flux (i.e., $c=q$, with $c_{\max}=q_{\max}$ and a nominal $\kappa=3$).

\begin{table}[hbt!]
\caption{\label{tab:entry-data} Data of the entry problem (angles in the aerodynamic and heating models in radians).}
\centering
\begin{tabular}{lll}
\hline\hline
Quantity & Symbol & Value \\
\hline
Initial altitude, velocity, flight-path angle & $h_0$, $v_0$, $\gamma_0$ & 260\,000 ft, 24\,061 ft/s, $-1.064$ deg \\
Initial longitude, latitude, heading & $\theta_0$, $\phi_0$, $\psi_0$ & 0, 0, 0 \\
Final altitude, velocity, flight-path angle & $h_f$, $v_f$, $\gamma_f$ & 80\,000 ft, 2\,500 ft/s, $-5$ deg \\
Final time & $t_f$ & 1\,500 to 2\,500 s \\
Bank angle, angle of attack & $\sigma$, $\alpha$ & $-90$ to 1 deg, $-90$ to 90 deg \\
Heating limit, nominal risk margin & $q_{\max}$, $\kappa$ & 70 BTU/(ft$^2$ s), 3 \\
\hline
Mass, reference area & $m$, $S$ & 7\,008 slug, 2\,690 ft$^2$ \\
Lift coefficients & $c_{L,0}$, $c_{L,1}$ & $-0.207$, 1.676 \\
Drag coefficients & $c_{D,0}$, $c_{D,1}$, $c_{D,2}$ & 0.07854, $-0.3529$, 2.04 \\
Common heating coefficient & $c_q$ & $9.289 \cdot 10^{-9}$ \\
Specific heating coefficients & $c_0$, $c_1$, $c_2$, $c_3$ & 1.067, $-1.101$, 0.6988, $-0.1903$ \\
\hline
Nominal density at zero altitude & $\bar{\rho}_0$ & $2.378 \cdot 10^{-3}$ slug/ft$^3$ \\
Nominal inverse scale height & $\bar{\beta}$ & $4.20168 \cdot 10^{-5}$ 1/ft \\
Earth radius, rotation rate & $R_{\oplus}$, $\omega_{\oplus}$ & 20\,902\,900 ft, $7.2722 \cdot 10^{-5}$ rad/s \\
Gravitational parameter & $\mu_{\oplus}$ & $1.4077 \cdot 10^{16}$ ft$^3$/s$^2$ \\
\hline\hline
\end{tabular}
\end{table}

\par
\begin{figure*}[htb]
    \centering
    \begin{subfigure}{3.25in}
        \centering
        \includegraphics[width=\linewidth]{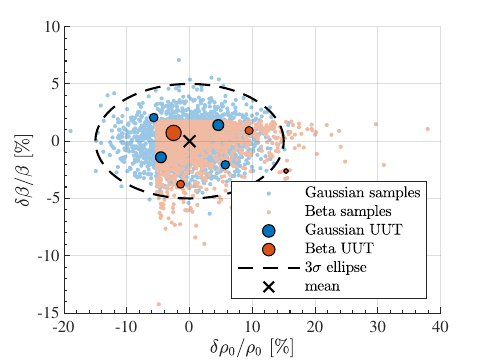}
        \caption{}
        \label{subfig:density-samples}
    \end{subfigure}%
    \hfill
    \begin{subfigure}{3.25in}
        \centering
        \includegraphics[width=\linewidth]{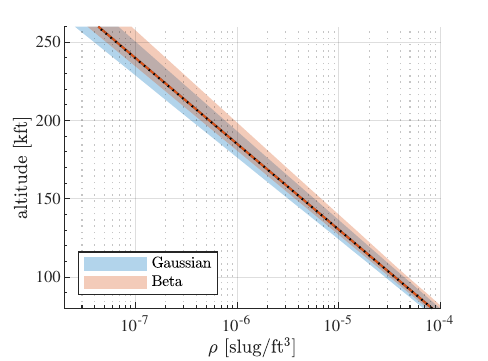}
        \caption{}
        \label{subfig:density-profiles}
    \end{subfigure}
    \caption{Gaussian and Beta density models with the same mean and covariance: a) samples and third-order rules of the G and B designs; b) density along the altitude.}
    \label{fig:density-models}
\end{figure*}

 The uncertain variables are the two parameters of Eq.~\eqref{eq:density-model}, $\mathbf{x} = [\rho_0, \beta]^T$: these are independent, with mean at the nominal values $\bar{\rho}_0$ and $\bar{\beta}$ of Table~\ref{tab:entry-data} and standard deviations equal to 5\% and (5/3)\% of these values (i.e., deviations of 15\% and 5\% at three standard deviations) \cite{schwartz2025_UnscentedHypersonicTrajectory}. To better showcase the relevance of high-order statistical moments, in this study we expand the uncertainty model of Ref.~\cite{schwartz2025_UnscentedHypersonicTrajectory} by considering two distributions that share this mean and covariance values. The Gaussian model is $\mathbf{x} \sim \mathcal{N}(\boldsymbol{\mu}, \mathbf{P})$, while in the Beta model each component is an affine map of a $\mathrm{Beta}(1, 20)$ variable, scaled to the same mean and variance, with the map of $\beta$ reversed such that both tails point towards a denser atmosphere. This model represents a somewhat "worst-case" scenario with respect to a conservative Gaussian assumption when only relying on mean and covariance information, as shown in Fig.~\ref{fig:density-models}.

\par The 2000 samples drawn from the two distributions, shown in Fig.~\ref{subfig:density-samples}, clearly show that despite sharing mean and covariance, the distributions represent very different shapes. The corresponding \ac{uut} rules are also shown, and we observe that the \ac{uut} places one point in the tail for the Beta distribution with a small weight, while all other points (for both the Gaussian and Beta rules) remain closer to the mean and with greater weights. The corresponding density profiles are shown in Fig.~\ref{subfig:density-profiles}: as expected, a large $\rho_0$ or a small $\beta$ raises the density at every altitude, therefore the Beta distribution, which is skewed towards those values, corresponds to a generally denser atmosphere than its Gaussian counterpart.

\par Five designs are compared, each named after the density model used in its guidance problem. The deterministic design D solves the problem with the nominal parameters only (i.e., Eq.~\eqref{eq:unscented-ocp} with $N = 1$, the point at the mean, and $\kappa = 0$). The unscented designs G and B solve Eq.~\eqref{eq:unscented-ocp} with \ac{uut} rules matching $\mathcal{J}_3^{\mathrm{f}}$ (i.e., third-order rules) of the Gaussian and Beta models, respectively, and the designs G4 and B4 with rules matching $\mathcal{J}_4^{\mathrm{f}}$ of the same models (i.e., fourth-order rules). Each rule is named after the design that uses it (e.g., the B4 rule). Table~\ref{tab:entry-designs} collects the five designs with their prescribed sets. The \ac{ug} results for each of these designs are then independently validated by \ac{mc} sampling of the density model, which does not necessarily match the one used to design the \ac{ug} solution: we denote by X/Y the design X validated on the model Y. The cases D/G and G/G, where the atmosphere is Gaussian, assess what the chance constraint formulation of Eq.~\eqref{eq:kappa-chance-constraint} gains over the nominal design (reproducing the setting of Ref.~\cite{schwartz2025_UnscentedHypersonicTrajectory}). The cases where the atmosphere follows the Beta model assess three further aspects: the effect of a Gaussian assumption when the actual distribution is skewed (G/B against B/B), whether matching more moments of the Gaussian model compensates for this assumption (G4/B), and what the \ac{uut} recovers by matching the actual moments up to the fourth order (B4/B). Note that the G/B against B/B comparison provides some insight into the potential advantages of adopting the \ac{uut} approach in uncertain \acp{ocp}: since we established in Sec.~\ref{sec:numerical-assessment} that existing sigma-point methods are particular solutions of the \ac{uut} problem, the G design under the B model corresponds to the adoption of an \ac{ut} rule of a beta distribution (i.e., neglecting its statistical moments except for mean and covariance).

\subsection{Results and Monte Carlo Validation}
\label{subsec:entry-results}

\par The \ac{ug} problems are transcribed and solved in \acs{raptor} \cite{minnozzi2026_TheoreticalComputationalAdvances}, using an hp-pseudospectral transcription on 32 Birkhoff \acl{flgr} segments of 10 nodes each: the problem is scaled using \ac{pjrn} and solved using \ac{ipopt} with a $10^{-8}$ tolerance. The chosen mesh is a trade-off between computational speed and node density, offering a highly accurate \ac{ps} representation of the solution at acceptable solution times. The \ac{mc} validation uses $5 \cdot 10^4$ samples for each density model, and the candidate optimal controls of each design are propagated using an adaptive Runge--Kutta integrator, with relative and absolute tolerances $10^{-8}$. The risk associated with each case (i.e., the combination of control/model) at threshold $\bar{q}$ is the percentage of sample trajectories whose peak heating rate exceeds $\bar{q}$. Since the optimal trajectories keep their heating profiles at exactly $q_{\max}$, computing the risk at exactly $\bar{q}=q_{\max}$ implies that most violations are caused by overshoots between collocation nodes, or minor differences caused by the independent numerical propagation: to exclude these effects from the analysis we follow the same approach proposed in Ref.~\cite{schwartz2025_UnscentedHypersonicTrajectory}, and consider risks at $\bar{q} = 72$~BTU/(ft$^2$ s), $3\%$ above the limit.

\par \par
\begin{figure*}[htb]
    \centering
    \begin{subfigure}{3.25in}
        \centering
        \includegraphics[width=\linewidth]{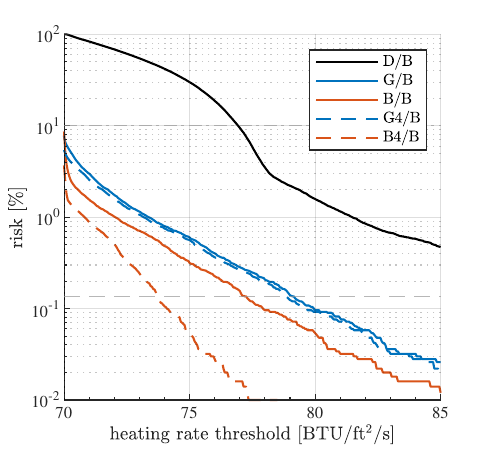}
        \caption{}
        \label{subfig:risk-beta}
    \end{subfigure}%
    \hfill
    \begin{subfigure}{3.25in}
        \centering
        \includegraphics[width=\linewidth]{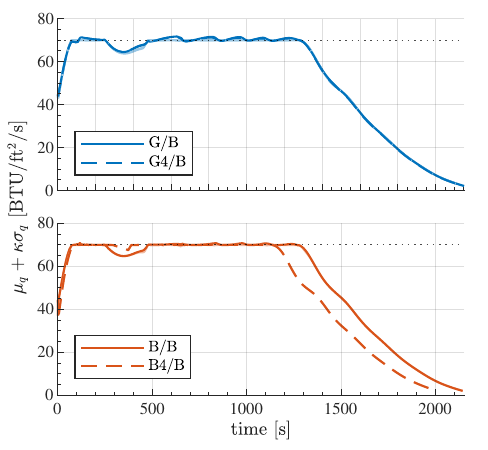}
        \caption{}
        \label{subfig:heating-beta}
    \end{subfigure}
    \caption{Designs with $\kappa = 3$ in the Beta atmosphere: a) risk chart; b) constrained statistic $\mu_q + 3\sigma_q$ measured by \ac{mc} sampling and predicted by each rule.}
    \label{fig:entry-nominal}
\end{figure*}

 The nominal results are computed for $\kappa=3$ using the ``deviation around the mean'' approach of Eq.~\eqref{eq:kappa-chance-constraint} and summarized in Table~\ref{tab:entry-designs}, reporting the size, objective, and computational cost of each design, while considerations on risk and the individual solutions are introduced through Fig.~\ref{fig:entry-nominal}. As expected, designs that account for uncertainties tend to reduce their performance in terms of crossrange (i.e., lower $|\phi(t_f)|$). The introduction of uncertainty also introduces a significant computational burden into the problem: each sigma point corresponds to a state trajectory, such that the number of states grows linearly with the number of points. All designs converge in a comparable number of \ac{ipopt} iterations, but the time per iteration grows faster than the number of states: with 4 and 6 times the states of the D design, the unscented designs take about 10 and 20 times longer per iteration, respectively.

\begin{table}[hbt!]
\caption{\label{tab:entry-designs} Entry designs with $\kappa = 3$: size, objective, and \ac{ipopt} iterations and solve time.}
\centering
\begin{tabular}{lcccccc}
\hline\hline
Design & Set & Points $N$ & States & $\phi(t_f)$, deg & Iterations & Solve time, s \\
\hline
D & -- & 1 & 6 & $-28.748$ & 684 & 75.4 \\
G & $\mathcal{J}_3^{\mathrm{f}}$ & 4 & 24 & $-28.053$ & 879 & 913 \\
B & $\mathcal{J}_3^{\mathrm{f}}$ & 4 & 24 & $-28.006$ & 890 & 926 \\
G4 & $\mathcal{J}_4^{\mathrm{f}}$ & 6 & 36 & $-28.049$ & 793 & 1713 \\
B4 & $\mathcal{J}_4^{\mathrm{f}}$ & 6 & 36 & $-27.724$ & 819 & 1752 \\
\hline\hline
\end{tabular}
\end{table}

\begin{figure*}[htb]
    \centering
    \begin{subfigure}{3.25in}
        \centering
        \includegraphics[width=\linewidth]{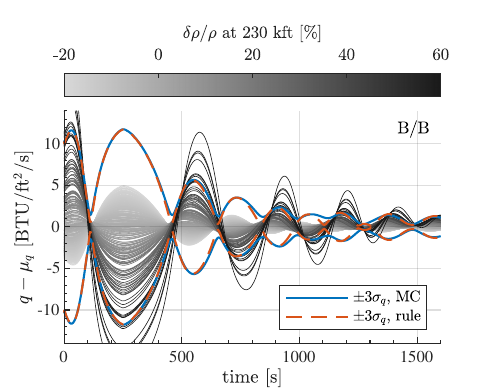}
        \caption{}
        \label{subfig:samples-bb}
    \end{subfigure}%
    \hfill
    \begin{subfigure}{3.25in}
        \centering
        \includegraphics[width=\linewidth]{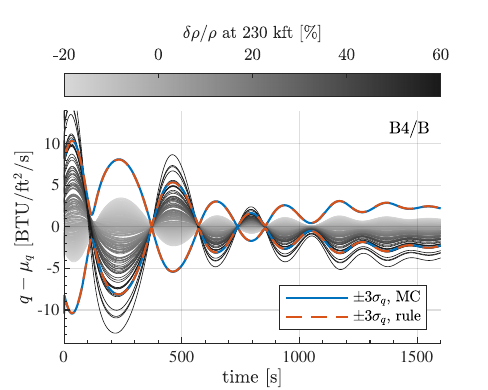}
        \caption{}
        \label{subfig:samples-b4b}
    \end{subfigure}
    \caption{Heating rate of 300 \ac{mc} samples about their mean in the Beta atmosphere, with the measured and predicted $\pm 3\sigma_q$ bands: a) B/B; b) B4/B.}
    \label{fig:heating-samples}
\end{figure*}

\par Figure~\ref{subfig:risk-beta} shows the risk chart of the designs in the Beta atmosphere: for each heating-rate threshold, the low-risk percentages of each design show its ability to address the uncertain atmospheric density. As shown by the blue lines, incorrectly assuming a Gaussian distribution produces higher-risk trajectories, and increasing the order of the \ac{uut} Gaussian approximation (i.e., adopting a higher-order scheme similar to the \ac{cut4}) yields negligible advantages. On the other hand, correctly matching the underlying distribution provides significant benefits: at $\bar{q}$, the risk of the B/B trajectories is about one third lower than that of the G4/B ones, and increasing the approximation order halves it again (i.e., 0.51\% for B4/B), as shown by the orange lines.

\par The main effect observed on the heating profiles (computed through Eq.~\eqref{eq:kappa-chance-constraint}) in Fig.~\ref{subfig:heating-beta} is the difference between the B4/B and B/B trajectories: the B solution has a dip between around 250--500~s and the final decrease around 1300~s, while the B4/B trajectory avoids the dip and starts decreasing earlier, around 1150~s. Figure~\ref{fig:heating-samples} explains this behavior through the heating rates of the \ac{mc} samples about their mean, with darker lines for denser atmospheres. The dip is caused by the common open-loop control: the samples in denser atmospheres generate more lift and climb above the others, such that the heating rates spread significantly after the first peak. Keeping the risk constraint at its bound therefore requires a lower mean heating rate, which lifts the mean trajectory: the statistic then remains below the bound, with the constraint inactive, until the descent brings the ensemble back to the bound around 480~s. Since the denser samples climb more (and receive less heating), their heating rates cross those of the less dense samples once the controls determine a descent arc (i.e., where dark and light lines swap in Fig.~\ref{fig:heating-samples}). At these crossings, the heating rate no longer varies linearly with the density, meaning that $\sigma_q$ is determined by the second-order terms, which only the fourth-order moments capture exactly (Sec.~\ref{subsec:moment-based-risk-constraint}). As shown by the bands in Fig.~\ref{fig:heating-samples}, the B rule underestimates $\sigma_q$ at each crossing, while the B4 rule, whose points reach further into the tail of the Beta distribution, matches the \ac{mc} value. The B/B design can therefore ride the bound until about 1300~s, incorrectly exceeding it at the crossing points, while the B4/B solution flies lower from the start with a smaller dispersion of the heating rate, avoiding the dip but leaving the bound earlier, sacrificing crossrange performance in favor of a safer trajectory.

\subsection{Effects of Higher-Order Information on Chance Constraints}
\label{subsec:skewness-margin}

\par A rule that correctly approximates $\sigma_q$ determines a low-risk trajectory, but it does not make the $\kappa=3$ value deliver the actual probability that the Gaussian assumption of Sec.~\ref{subsec:moment-based-risk-constraint} associates with it. Figure~\ref{fig:entry-nominal} shows that for $\kappa=3$, the trajectory risk at $\bar{q}=72$BTU/(ft$^2$ s) is $0.51\%$, against the target $0.135\%$ for the Gaussian-distributed constraint: the heating rate is not Gaussian, and a margin built only on its mean and standard deviation cannot account for its shape. We therefore compare the two constraints of Sec.~\ref{subsec:moment-based-risk-constraint} at the same target probability: the "deviation around the mean" constraint of Eq.~\eqref{eq:kappa-chance-constraint} with $\kappa = z$, and the skewness-based constraint of Eq.~\eqref{eq:cf-constraint}. For each combination of design/model (e.g., G/B, B4/B, and so on), two problems are solved, one with each formulation of the heating rate chance constraint, using the same \acs{raptor} settings presented previously. The chosen $z$ values are $z = 2$, $2.5$, $3$, and $3.5$, corresponding to acceptable violation probabilities of $\varepsilon = 2.28$, $0.621$, $0.135$, and $0.0233\%$. To better capture the small violation probabilities, the \ac{mc} validation for this analysis uses $2 \cdot 10^5$ samples for each design/model combination. The validation measures two quantities: the largest pointwise exceedance probability of $q_{\max}$ (i.e., the quantity that both constraints bound), and the risk at $\bar{q} = 72$~BTU/(ft$^2$ s), which counts the samples whose heating rate exceeds $\bar{q}$ at any time along the trajectory.

\begin{figure*}[htb]
    \centering
    \begin{subfigure}{3.25in}
        \centering
        \includegraphics[width=\linewidth]{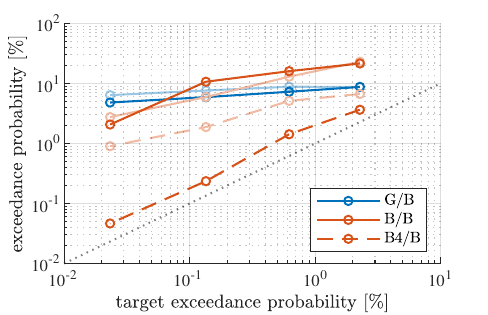}
        \caption{}
        \label{subfig:cf-calibration}
    \end{subfigure}%
    \hfill
    \begin{subfigure}{3.25in}
        \centering
        \includegraphics[width=\linewidth]{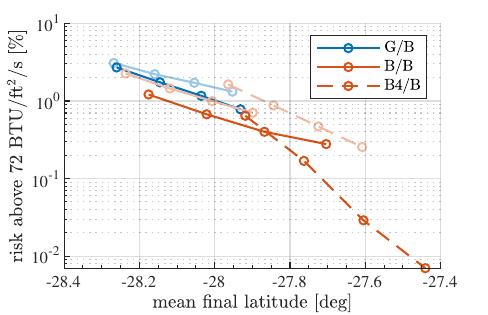}
        \caption{}
        \label{subfig:cf-tradeoff}
    \end{subfigure}
    \caption{Covariance-only and skewness-based margins in the Beta atmosphere: a) measured against target exceedance probability; b) risk at 72~BTU/(ft$^2$ s) against mean final latitude.}
    \label{fig:entry-skewness-margin}
\end{figure*}

\par Figure~\ref{subfig:cf-calibration} compares the measured and target probabilities, such that an exact link lies on the dotted diagonal. Each design appears as a pair of curves with the same color and line style, with semi-transparent lines for the covariance-only margin and opaque lines for the skewness-based one. Two distinct effects can be observed in the figure: the more accurate statistical information built from the 4th-order moments using the \ac{uut} is shown by the B4 risks being consistently lower than B (regardless of the inclusion or not of the third-order statistical information in the risk constraint), while the more precise risk approximation guaranteed by the third-order information is evidenced by the opaque B4 points being significantly closer to the target trend represented by the dotted diagonal. Noticeably, the gains caused by introducing higher-order information into the risk constraint are not reflected on the B and G designs, as in these cases the dominant source of error is the limited propagation accuracy, rather than the low-order approximation of the risk.

\par As one would expect, safer trajectories are not free: increasing the allowed exceedance probability (i.e., riskier trajectories) allows the optimization to improve the crossrange performance through riskier trajectories. Figure~\ref{subfig:cf-tradeoff} compares the designs in terms of trajectory risk against the achieved objective, through the same pairs of curves as Fig.~\ref{subfig:cf-calibration}, with $z$ increasing from left to right. The curves represent Pareto fronts for the competing problem objective (crossrange maximization) with the design safety constraint. For every design, the skewness-based curve lies below the covariance-only one, meaning that at the same crossrange the skewness-based margin yields a lower risk. The B4/B curves are also much steeper than other designs: since the B4 rule predicts the statistics of the heating rate accurately, a larger margin acts directly on the actual tail of its distribution, while for the other designs part of the exceedance comes from errors that the margin does not account for. At large crossrange, the B/B designs reach the lowest risk above 72~BTU/(ft$^2$ s): this is a consequence of the chosen threshold, as these trajectories frequently exceed the "original" $q_{\max}$ but by small amounts (e.g., 60\% of the samples exceed 70~BTU/(ft$^2$ s) with $z = 2$, against 4.4\% for B4/B).

\section{Conclusions}
\label{sec:conclusions}

\par This paper introduced the \ac{uut}, a sigma-point construction method that jointly optimizes point locations and weights, subject to prescribed moment-matching conditions, nonnegative weights with unit sum, and the support of the target distribution. Instead of fixing a point structure and deriving the matched moments from it, the \ac{uut} only requires the statistical moments of the target distribution and its support set, while the number of points, their locations and weights are chosen to provide a valid rule that is sparse (i.e., contains few active sigma-points) and compact (i.e., concentrates its points close to the mean).

\par The numerical assessment showed that the existing rules are particular solutions of the \ac{uut} problem wherever they are verified, while the \ac{uut} produces verified rules on most of the targets where existing methods fail. The proposed sparsification, which combines reweighting with fixed points and greedy elimination with free points, matches or outperforms a direct search over the number of points while producing compact rules. Evaluating the technique over multiple distributions and increasing dimensions showed that the \ac{uut} can match the moments of the \ac{ut}, \ac{cut4}, and \ac{gut} with fewer points, or more moments with the same number of points.

\par In the entry problem under atmospheric uncertainty, it is shown that a Gaussian assumption on the underlying distribution can produce higher-risk trajectories than using a more accurate \ac{uut} sigma-point scheme, and neither a larger constraint margin nor more Gaussian moments close this gap. The \ac{mc} validation traced the gap to the statistics that the rule reconstructs: only the rule matching fourth-order moments of the underlying beta distribution reproduces the standard deviation and skewness of the heating rate chance constraint. The problem also showcases how the higher-order information available from such construction can be used to improve the chance constraint prediction, bringing the measured violation probability within a factor two of its target value, with up to an order of magnitude less risk at the same crossrange performance.

\section*{Appendix}
\label{sec:appendix}

\setcounter{subsection}{0}
\renewcommand{\thesubsection}{\Alph{subsection}}

\subsection{Analytic Derivatives of the Moment Conditions}
\label{subsec:analytic-derivatives}

\par \ac{ipopt} receives the exact first derivatives of every constraint and of the objective of Eq.~\eqref{eq:rule-design-problem} with respect to the point locations and weights. The unit mass constraint has derivative one with respect to each weight and zero with respect to the points, and the support constraints of Eq.~\eqref{eq:whitened-support} are linear in the whitened points, with the rows of $\mathbf{L}$ as derivatives. When whitening is active, the decision variables are the whitened points $\mathbf{z}_i$, and the formulas below apply with $\mathbf{z}_i$ in place of $\mathbf{x}_i$ and with the whitened entries $m^{z}_{\mathbf{j}}$ of Eq.~\eqref{eq:whitened-moments} as targets. Only the moment conditions and the compactness objective of Eq.~\eqref{eq:compactness-objective} therefore require a derivation.

\par The first-order entries are the components of the mean $\hat{\boldsymbol{\mu}} = \mathbf{X}\mathbf{w}$, with derivatives
\begin{equation}
    \frac{\partial \hat{\mu}_p}{\partial w_i} = x_{i,p}, \qquad \frac{\partial \hat{\mu}_p}{\partial x_{i,q}} = w_i\, \delta_{pq}
    \label{eq:mean-derivatives}
\end{equation}
where $\delta_{pq}$ is the Kronecker delta.

\par For an entry of order $k \geq 2$, we introduce the centered product of sigma-point $i$
\begin{equation}
    e_{i,\mathbf{j}} = \prod_{l=1}^{k} \left(x_{i,j_l} - \hat{\mu}_{j_l}\right)
    \label{eq:centered-product}
\end{equation}
such that  Eq.~\eqref{eq:rule-moments} becomes $\hat{m}_{\mathbf{j}} = \sum_i w_i e_{i,\mathbf{j}}$. We also define the index tuple obtained by removing the $l$-th index from $\mathbf{j}$ as $\mathbf{j}_{\setminus l}$, and the weighted sum of the corresponding centered products as $s_{\mathbf{j}_{\setminus l}} = \sum_{r} w_r\, e_{r,\mathbf{j}_{\setminus l}}$. Differentiating each factor of Eq.~\eqref{eq:centered-product}, including its dependence on $\hat{\boldsymbol{\mu}}$ through Eq.~\eqref{eq:mean-derivatives}, yields
\begin{equation}
    \frac{\partial \hat{m}_{\mathbf{j}}}{\partial w_i} = e_{i,\mathbf{j}} - \sum_{l=1}^{k} s_{\mathbf{j}_{\setminus l}}\, x_{i,j_l}, \qquad
    \frac{\partial \hat{m}_{\mathbf{j}}}{\partial x_{i,p}} = w_i \sum_{l \,:\, j_l = p} \left(e_{i,\mathbf{j}_{\setminus l}} - s_{\mathbf{j}_{\setminus l}}\right)
    \label{eq:moment-derivatives}
\end{equation}
where the last sum runs over the positions $l$ of $\mathbf{j}$ whose index equals $p$. In both derivatives, the first term is the direct effect of the weight or point location on the centered product for point $i$, while the terms with $s_{\mathbf{j}_{\setminus l}}$ account for the shift in the rule that the same weight or location produces. For instance, for $\mathbf{j} = (1,1,2)$, the tuples $\mathbf{j}_{\setminus 1}$ and $\mathbf{j}_{\setminus 2}$ are both $(1,2)$, and $\mathbf{j}_{\setminus 3} = (1,1)$, such that
\begin{equation}
    \frac{\partial \hat{m}_{(1,1,2)}}{\partial x_{i,1}} = 2 w_i \left(e_{i,(1,2)} - s_{(1,2)}\right), \qquad \frac{\partial \hat{m}_{(1,1,2)}}{\partial x_{i,2}} = w_i \left(e_{i,(1,1)} - s_{(1,1)}\right)
    \label{eq:moment-derivatives-example}
\end{equation}
where the factor two comes from the repeated index. For $k \geq 3$, each sum $s_{\mathbf{j}_{\setminus l}}$ is the rule entry $\hat{m}_{\mathbf{j}_{\setminus l}}$ of order $k-1$, independently from its presence. For $k = 2$, it is the weighted sum of the centered coordinates, which goes to zero as the unit-mass constraint becomes satisfied.

\par The compactness objective of Eq.~\eqref{eq:compactness-objective} is a weighted sum of powers of the whitened radii, whose derivatives with respect to the whitened points and the weights are
\begin{equation}
    \frac{\partial J_{\mathrm{c}}}{\partial w_i} = \frac{\left\|\mathbf{z}_i\right\|^{q}}{D}, \qquad \frac{\partial J_{\mathrm{c}}}{\partial \mathbf{z}_i} = \frac{q\, w_i \left\|\mathbf{z}_i\right\|^{q-2}}{D}\, \mathbf{z}_i
    \label{eq:compactness-derivatives}
\end{equation}
where $D = \sum_i \bar{w}_i \|\bar{\mathbf{z}}_i\|^{q}$ is the constant denominator of Eq.~\eqref{eq:compactness-objective}.

\subsection{Effect of Plain Weight Regularization}
\label{subsec:regularization-evidence}

\par
\begin{figure*}[htb]
    \centering
    \begin{subfigure}{3.25in}
        \centering
        \includegraphics[width=\linewidth]{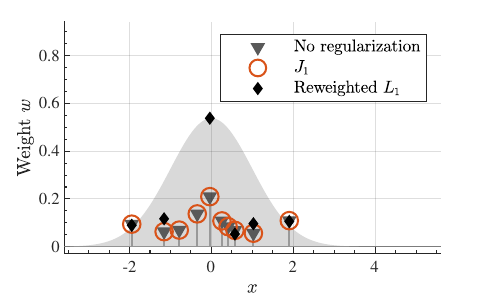}
        \caption{}
        \label{subfig:regularization-l1}
    \end{subfigure}%
    \hfill
    \begin{subfigure}{3.25in}
        \centering
        \includegraphics[width=\linewidth]{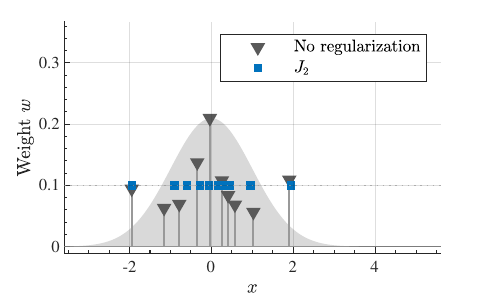}
        \caption{}
        \label{subfig:regularization-l2}
    \end{subfigure}
    \caption{Plain weight regularization of a ten-point rule for the scalar Gaussian target: a) $J_1$ and reweighted $L_1$; b) $J_2$.}
    \label{fig:regularization-evidence}
\end{figure*}

 Section~\ref{subsec:reweighted-l1} opens with the argument that plain weight regularization does not provide the sparsity-inducing capability required for the \ac{uut} framework. A numerical showcase of the effects of plain $L_1$ and $L_2$ regularization is provided in Fig.~\ref{fig:regularization-evidence}, comparing the sparsification obtained by the reweighted $L_1$ approach illustrated in Alg.~\ref{alg:reweighting} (i.e., the first step of the \ac{uut} sparsification pipeline) against the plain nonlinear solves using $J_1$ and $J_2$ terms of Eq.~\eqref{eq:plain-regularization}. The chosen target is a scalar Gaussian distribution, enforcing moment matching on the $\mathcal{J}_5^{\mathrm{f}}$ set starting from an overparameterized rule with $N=10$ points.

\par Figure~\ref{subfig:regularization-l1} exhibits the foreseen behavior for the $L_1$ terms: $J_1$ cannot promote sparsity effectively, as the unit-mass constraint prevents it from having any effects, while the reweighted $L_1$, despite the fixed-points constraints, successfully pushes redundant weights to 0, identifying a sigma-point rule exactly on its theoretical $N_C$ bound. The $J_2$ effect is clear in Fig.~\ref{subfig:regularization-l2}, as it promotes a uniformly weighted rule (with weights $w_i=1/N$) with points clustered close to the mean. While this property cannot effectively push for the sparsity desired in applications of the \ac{uut} approach, if the sigma-point budget $N$ is prescribed or reliably estimated a-priori, the push for an equal-weighted rule is a desirable feature for \ac{ug} problems, as it generates well conditioned statistical reconstruction among the propagated samples (i.e., the problem's sensitivity to each uncertainty sample is spread evenly). Equal weights, however, are reached only when the moment conditions allow them, which is not always the case (e.g., skewed targets). Since it does not reduce the number of points and conflicts with the skewed targets, the $L_2$ term is not used in Sec.~\ref{sec:numerical-assessment}, where the compactness objective of Eq.~\eqref{eq:compactness-objective} keeps the sparsified rules away from the tails instead.

\subsection{Moments of the Target Distributions}
\label{subsec:target-moments}

\par The moment conditions of Eq.~\eqref{eq:moment-conditions} require the central moments $m_{\mathbf{j}}$ of Eq.~\eqref{eq:target-moments} of each target of Table~\ref{tab:targets}, up to order eight. All of them are available in closed form or as finite sums, such that neither sampling nor numerical quadrature enters the comparisons of Sec.~\ref{sec:numerical-assessment}. This appendix collects them in one place, starting from the scalar laws and then combining them into the multivariate entries.

\par Table~\ref{tab:target-moments} collects, for each scalar law, its support, its moments $r_q = \mathbb{E}[x^q]$ of any order $q$ (on the second line of each law), and the four statistics that summarize its shape: mean $\mu$, variance $\sigma^2$, skewness, and kurtosis. In the table, $(q-1)!! = 1 \cdot 3 \cdots (q-1)$, and $r_q$ vanishes for odd $q$ for the two symmetric laws. Each target is shown in standardized form $z = (x - \mu)/\sigma$, with the first four moments being therefore 0, 1, and the skewness and kurtosis of Table~\ref{tab:target-moments} \cite{johnson1994_ContinuousUnivariateDistributions}. The moments of higher order follow from $r_q$ by expanding $(x - \mu)^k$ with the binomial theorem,
\begin{equation}
    \mathbb{E}\left[z^k\right] = \frac{1}{\sigma^k} \sum_{q=0}^{k} \binom{k}{q} (-\mu)^{k-q} r_q
    \label{eq:standardized-moments}
\end{equation}
which reduces to $\mathbb{E}[z^k] = r_k/\sigma^k$ for the laws with zero mean.

\begin{table}[hbt!]
\caption{\label{tab:target-moments} Support, mean, variance, skewness, kurtosis, and moments $r_q$ of any order (second line) of the scalar laws of Table~\ref{tab:targets}.}
\centering
\setlength{\tabcolsep}{4pt}
\renewcommand{\arraystretch}{2.2}
\begin{tabular}{>{$}l<{$}>{$}l<{$}>{$\displaystyle}l<{$}>{$\displaystyle}l<{$}>{$\displaystyle}l<{$}>{$\displaystyle}l<{$}}
\hline\hline
\text{Law} & \text{Support} & \mu & \sigma^2 & \mathbb{E}[(x - \mu)^3]/\sigma^3 & \mathbb{E}[(x - \mu)^4]/\sigma^4 \\
\hline
\mathcal{N}(0, 1) & \mathbb{R} & 0 & 1 & 0 & 3 \\
 & \multicolumn{5}{>{$\displaystyle}l<{$}}{r_q = (q-1)!! \text{ for even } q} \\
\mathcal{U}(0, 1) & [0, 1] & \frac{1}{2} & \frac{1}{12} & 0 & \frac{9}{5} \\
 & \multicolumn{5}{>{$\displaystyle}l<{$}}{r_q = \frac{1}{q+1}} \\
\mathrm{Beta}(2, b) & [0, 1] & \frac{2}{b+2} & \frac{2b}{(b+2)^2(b+3)} & \frac{2(b-2)\sqrt{b+3}}{(b+4)\sqrt{2b}} & 3 + \frac{3\left[(b-2)^2(b+3) - 2b(b+4)\right]}{b(b+4)(b+5)} \\
 & \multicolumn{5}{>{$\displaystyle}l<{$}}{r_q = \prod_{l=0}^{q-1} \frac{2+l}{b+2+l}} \\
\mathrm{Gamma}(\alpha, 1) & [0, \infty) & \alpha & \alpha & \frac{2}{\sqrt{\alpha}} & 3 + \frac{6}{\alpha} \\
 & \multicolumn{5}{>{$\displaystyle}l<{$}}{r_q = \frac{\Gamma(\alpha+q)}{\Gamma(\alpha)}} \\
\mathrm{LogN}(0, s^2) & (0, \infty) & e^{s^2/2} & \left(e^{s^2} - 1\right)e^{s^2} & \left(e^{s^2} + 2\right)\sqrt{e^{s^2} - 1} & e^{4s^2} + 2e^{3s^2} + 3e^{2s^2} - 3 \\
 & \multicolumn{5}{>{$\displaystyle}l<{$}}{r_q = e^{q^2 s^2/2}} \\
\mathrm{Weibull}(c, 1) & [0, \infty) & \Gamma_1 & \Gamma_2 - \Gamma_1^2 & \frac{\Gamma_3 - 3\Gamma_1\Gamma_2 + 2\Gamma_1^3}{\sigma^3} & \frac{\Gamma_4 - 4\Gamma_1\Gamma_3 + 6\Gamma_1^2\Gamma_2 - 3\Gamma_1^4}{\sigma^4} \\
 & \multicolumn{5}{>{$\displaystyle}l<{$}}{r_q = \Gamma_q = \Gamma\!\left(1 + \frac{q}{c}\right)} \\
\frac{1}{2}\mathcal{N}\big({\pm}d, \frac{1}{4}\big) & \mathbb{R} & 0 & d^2 + \frac{1}{4} & 0 & \frac{d^4 + \frac{3}{2}d^2 + \frac{3}{16}}{\left(d^2 + \frac{1}{4}\right)^2} \\
 & \multicolumn{5}{>{$\displaystyle}l<{$}}{r_q = \sum_{l \text{ even}} \binom{q}{l} \frac{d^{q-l}\,(l-1)!!}{2^{l}} \text{ for even } q} \\
\hline\hline
\end{tabular}
\end{table}

\par The multivariate entries follow from the scalar moments in one of three ways, depending on the relation between the components of the target. When the components are independent, the expectation of each product in Eq.~\eqref{eq:target-moments} factors into the product of the expectations of its components,
\begin{equation}
    m_{\mathbf{j}} = \prod_{i=1}^{n} \mathbb{E}\left[z^{c_i(\mathbf{j})}\right]
    \label{eq:independent-moments}
\end{equation}
where $c_i(\mathbf{j})$ is the number of times the index $i$ appears in $\mathbf{j}$, and $\mathbb{E}[z^0] = 1$. For instance, $m_{(1,1,2,2)} = \mathbb{E}[z^2]^2 = 1$, while $m_{(1,1,2)} = \mathbb{E}[z^2]\,\mathbb{E}[z] = 0$: every entry in which a component appears exactly once is therefore zero. For the correlated Gaussian targets, the entries of odd order vanish, while each entry of even order is the sum, over all the ways of grouping its $k$ indices into pairs, of the products of the covariance entries of the pairs \cite{isserlis1918_FormulaProductMomentCoefficient}, computed as
\begin{equation}
    m_{\mathbf{j}} = \sum_{\pi \in \Pi_k} \prod_{(a, b) \in \pi} \left[\mathbf{R}_{\rho}\right]_{j_a j_b}
    \label{eq:isserlis}
\end{equation}
where $\Pi_k$ is the set of the $(k-1)!!$ groupings of the positions $\{1, \dots, k\}$ into pairs. For instance, $m_{(1,1,2,2)} = 1 + 2\rho^2$, from the pairings $\{(1,2),(3,4)\}$, $\{(1,3),(2,4)\}$, and $\{(1,4),(2,3)\}$. Finally, the mixed exponentials $\mathbf{x} = \mathbf{R}_{\rho}^{1/2}\boldsymbol{\epsilon}$ are a linear map of the vector $\boldsymbol{\epsilon}$ of independent standardized exponentials, whose entries $m^{\epsilon}_{\mathbf{q}}$ follow from Eq.~\eqref{eq:independent-moments}. Each entry of $\mathbf{x}$ is therefore a weighted sum of the entries of $\boldsymbol{\epsilon}$ of the same order,
\begin{equation}
    m_{\mathbf{j}} = \sum_{\mathbf{q}} \left(\prod_{l=1}^{k} \left[\mathbf{R}_{\rho}^{1/2}\right]_{j_l q_l}\right) m^{\epsilon}_{\mathbf{q}}
    \label{eq:mixed-moments}
\end{equation}
as for the whitening of Eq.~\eqref{eq:whitened-moments}, where the sum runs over all the index tuples $\mathbf{q}$ of length $k$. For instance, since the only nonzero third-order entries of $\boldsymbol{\epsilon}$ are its marginal skewnesses, equal to 2, the mixed entry $m_{(1,1,2)} = 2\sum_{i=1}^{n} [\mathbf{R}_{\rho}^{1/2}]_{1i}^2 [\mathbf{R}_{\rho}^{1/2}]_{2i}$ is nonzero, and equals 0.61 for $n = 2$.

\subsection{Deterministic Initial Guess}
\label{subsec:initial-guess}

\par The initial guess of Sec.~\ref{subsec:rule-design-problem} places $N$ equally weighted points on a Gaussian surrogate of the target, with the same mean and covariance, and it is built in three steps. First, $N$ points $\mathbf{v}_i \in (0,1)^n$ of a Sobol sequence \cite{sobol1967_DistributionPointsCube} are mapped to standard Gaussian samples
\begin{equation}
    \tilde{\mathbf{z}}_i = \Phi^{-1}(\mathbf{v}_i), \qquad i = 1, \dots, N
    \label{eq:sobol-gaussian}
\end{equation}
where $\Phi^{-1}$ is the inverse cumulative distribution function of the standard Gaussian distribution. The first points of the Sobol sequence are skipped, since $\Phi^{-1}(0)$ is infinite. Second, the samples are centered on their mean and multiplied by the inverse Cholesky factor of their covariance (both computed according to Eq.~\eqref{eq:sigma-point-rule} with the equal weights $1/N$), which yields points $\mathbf{z}_i$ with exactly zero mean and identity covariance. These points are the initial guess in the whitened variables of Eq.~\eqref{eq:whitening}. Third, the points are mapped back to the original variables through
\begin{equation}
    \mathbf{x}_i = \boldsymbol{\mu} + \mathbf{L}\mathbf{z}_i, \qquad i = 1, \dots, N
    \label{eq:initial-guess}
\end{equation}
where $\mathbf{L}$ is the Cholesky factor of the target covariance of Sec.~\ref{subsec:problem-scaling}, such that the initial rule matches the mean and covariance of the target exactly. This construction avoids random sampling, improving the consistency and reproducibility of the results over the more "immediate" approach of using randomly sampled points from a Gaussian distribution. The retries of the greedy elimination of Sec.~\ref{subsec:greedy-elimination} perturb the points along directions taken from later points of the same sequence: each retry generates samples from the Sobol sequence as outlined above and perturbs the points of the starting rule by a factor 0.1 along the directions of the generated samples.

\section*{Acknowledgments}
\label{sec:acknowledgments}
\par AI was used to proofread and polish the manuscript text and formatting of the results. All content was reviewed, verified, and revised by the authors, who take full responsibility for the paper.
\par Funded by the European Union (\ac{erc}, STARGATE, 101169959). Views and opinions expressed are however those of the author(s) only and do not necessarily reflect those of the European Union or the European Research Council. Neither the European Union nor the granting authority can be held responsible for them.
\begin{figure}[h!]
\centering
\includegraphics[width=3.25in]{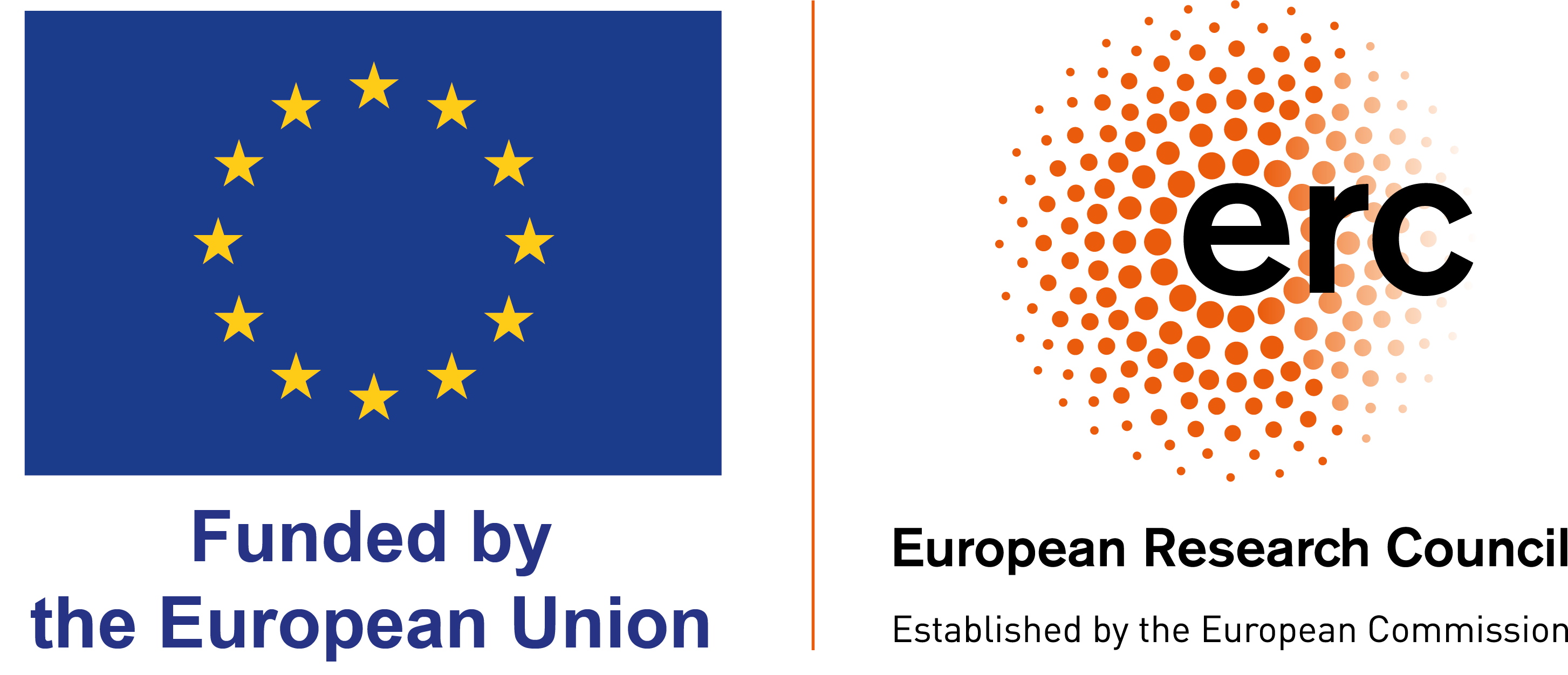}
\captionsetup{labelformat=empty}
\caption{}
\end{figure}

\bibliography{references.bib}

\end{document}